\documentclass[a4paper,twoside,11pt]{article}

\usepackage[utf8]{inputenc}   
\usepackage{lmodern}
\usepackage[T1]{fontenc}     
\usepackage[francais]{babel}  
\usepackage[a4paper]{geometry}

\usepackage{amssymb}
\usepackage{amsmath}
\usepackage{graphics}
\usepackage{epsfig}
\usepackage{graphicx}
\usepackage{lscape}

\usepackage{encadrement}
\usepackage[nottoc, notlof, notlot]{tocbibind}

\newboxedtheorem[boxcolor=black, background=white, titlebackground=blue!20,titleboxcolor = black]{res}{Résumé}{}
\newboxedtheorem[boxcolor=red, background=blue!5, titlebackground=blue!20,titleboxcolor = black]{theo}{Théorème}{}

\newboxedtheorem[boxcolor=orange, background=blue!3, titlebackground=blue!15,titleboxcolor = black]{prop}{Proposition}{}

\newboxedtheorem[boxcolor=green, background=blue!1, titlebackground=blue!10,titleboxcolor = black]{lemme}{Lemme}{lemme}

\newboxedtheorem[boxcolor=black, background=blue!2, titlebackground=blue!10,titleboxcolor = black]{methode}{Principe de la méthode}{}

\newboxedtheorem[boxcolor=black, background=blue!2, titlebackground=blue!10,titleboxcolor = black]{retenir}{Retenons!!}{}

\newboxedtheorem[boxcolor=black, background=blue!2, titlebackground=blue!10,titleboxcolor = black]{defin}{Définition}{}

\newboxedtheorem[boxcolor=red, background=blue!2, titlebackground=blue!10,titleboxcolor = black]{corollaire}{Corollaire}{}

\newboxedtheorem[boxcolor=black, background=blue!2, titlebackground=blue!5,titleboxcolor = black]{remarque}{Remarque}{}

\DeclareMathOperator{\Id}{Id} 

\usepackage{url}
\usepackage{fancyhdr}
\title{\textsc{Autour de la décomposition de Dunford réelle ou complexe :} \\
\textit{Théorie spectrale et méthodes effectives.}}     

\author{\textsc{Alaeddine BEN RHOUMA}\\
Professeur agrégé de Mathématiques \\
\\
\url{alaeddine.ben-rhouma@ac-guyane.fr}}

\date{ Juillet 2013}

\begin{document}

\maketitle


\begin{res}
Ces notes ne prétendent pas se substituer à un cours d'algèbre linéaire sur la réduction des endomorphismes ni à une présentation exhaustive de la décomposition de Dunford.

On se limitera au cas où le corps de base est $\mathbb{R}$ ou $\mathbb{C}$, et l'objectif de cette présentation est de faire un état des lieux sur les diverses méthodes de décomposition de Dunford.

 Lorsque les valeurs propres sont connues avec leurs valeurs exactes, une décomposition en éléments simples de l'inverse d'un polynôme annulateur nous fournit les projecteurs spectraux et a fortiori la décomposition attendue.

Le cas le plus délicat se présente dès que le spectre de l'endomorphisme n'est pas à notre disposition ce qui est une situation courante lorsque la dimension de l'espace vectoriel est supérieure à 4.

La méthode de Newton-Raphson vient alors à la rescousse pour nous fournir une suite qui converge quadratiquement vers la composante diagonalisable.

Certes, cette méthode très en vogue est assez efficace quelle que soit la taille de la matrice étudiée, mais elle nous laisse sur notre faim. En effet, on sait que les composantes de Dunford sont des polynômes en la matrice et on souhaiterait connaître ces polynômes générateurs.

La bonne nouvelle, c'est qu'une méthode effective utilisant le lemme chinois existe et elle a été introduite par Chevalley dans les années cinquantes du siècle dernier. 

Je mettrai l'accent sur cette méthode qui a été évoquée dans un article de Danielle Couty, Jean Esterle et Rachid Zarouf, en détaillant la preuve de l'algorithme dans le cas où le polynôme caractéristique est scindé sur le corps de base, puis je détaillerai le cas réel qui est une situation plus subtile nécéssitant  une étude plus poussée.

Un rappel sur les endomorphismes semi-simples a été introduit afin de justifier l'importance de la recherche d'une méthode effective permettant de tester la diagonalisabilité dans $\mathcal{M}_n(\mathbb{R})$ lorsqu'on ne dispose pas des valeurs propres de l'endomorphisme étudié. Pour y parvenir j'ai proposé les suites de Sturm comme outil de vérification de diagonalisabilité dans $\mathbb{R}$.

\end{res}


\newpage

\tableofcontents              

\newpage

\section{Résultats préliminaires}              

$\mathbb{K}$ désigne le corps $\mathbb{R}$ ou $\mathbb{C}$ et $E$ un espace vectoriel sur $\mathbb{K}$ de dimension $n\in\mathbb{N}^*$.

Pour tout $u\in \mathcal{L}(E)$, on note $\chi_u=\det (u-X \Id)$ le polynôme caractéristique de $u$ , $\mu_u$ son polynôme minimal et $\mathbb{K}[u]=\{P(u) ; P\in\mathbb{K}[X]\}$ qui est une algèbre commutative.

\subsection{Lemme des noyaux et ses conséquences} 
\subsubsection{Lemme des noyaux} 
 \begin{theo}[Décomposition des noyaux]

Soit $u\in\mathcal{L}(E)$ et $P=P_1...P_k\in\mathbb{K}[X]$ tel que $P_1,...,P_k$ soient premiers entre eux deux à deux.

Alors : $\ker P(u) = \bigoplus\limits_{i=1}^k \ker P_i(u)$.

\end{theo}

\paragraph{Preuve:}

On procède par récurrence sur $k\in\mathbb{N}$ et $k\geq2$.
\\
\\
\underline{Initialisation:} Soit $P=P_1P_2\in\mathbb{K}[X]$ tel que $P_1\wedge P_2=1$. D'après le théorème de Bezout, il existe $U,V\in\mathbb{K}[X]$ tels que $UP_1+VP_2=1$.

 On a alors pour tout $x\in E$ :

$$x=U(u)\circ P_1(u)(x)+V(u)\circ P_2(u)(x)$$

Soit $x\in\ker P_1(u)\cap \ker P_2(u)$. On obtient dans ce cas :

$u(u)\circ P_1(u)(x) = V(u)\circ P_2(u)(x) = 0$ et donc $x = 0$.

On a alors $\ker P_1(u)\cap \ker P_2(u) = \{0\}$.

Or si $x\in\ker P(u)$, avec $x=U(u)\circ P_1(u)(x)+V(u)\circ P_2(u)(x)$, 

on voit que $P_2(u)\Bigl(U(u)\circ P_1(u)(x)\Bigr)=0$ et $P_1(u)\Bigl(V(u)\circ P_2(u)(x)\Bigr)=0$.

On a alors, $x=x_1+x_2$ avec $x_1\in\ker P_2(u)$ et $x_2\in\ker P_1(u)$.

Donc $\ker P(u) = \ker P_1(u) \bigoplus \ker P_2(u)$.
\\
\\
\underline{Hérédité:} Supposons la propriété vraie jusqu'au rang $k$ et montrons qu'elle l'est encore au rang $k+1$.

Soit alors $P_1, \ldots ,P_k,P_{k+1}\in\mathbb{K}[X]$ premiers entre eux deux à deux.

Soit $P=P_1\cdots P_kP_{k+1}$. Puisque $P_1\cdots P_k$ et $P_{k+1}$ sont premiers entre eux, d'après l'intialisation on a $\ker P(u)=\ker (P_1\cdots P_k)(u)\bigoplus\ker P_{k+1}(u)$. 

D'après l'hypothèse de récurrence $\ker (P_1\cdots P_k)(u)=\bigoplus\limits_{i=1}^k\ker P_i$.

 Donc $\ker P(u)=\bigoplus\limits_{i=1}^{k+1}\ker P_i(u)$. $\blacksquare$

\subsubsection{Décomposition de l'espace en somme de sous-espaces caractéristiques} 

 \begin{prop}[Décomposition spectrale]

Si un élément  $P$ de $\mathbb{K}[X]$ est un polynôme annulateur d'un élément $u$ de $\mathcal{L}(E)$ tel que : $P(X)=\prod\limits_{i=1}^r (X-\lambda_i)^{\alpha_i}$, 

avec pour tout $1\leq i,j \leq r$, $\lambda_i \in \mathbb{K}$  , $\alpha_i\geq1$ et $\lambda_i \neq \lambda_j$ pour $i\neq j$. 

Alors $E=\bigoplus\limits_{i=1}^r \ker(u-\lambda_i \Id)^{\alpha_i}$.

\end{prop}

\paragraph{Preuve:}

C'est une conséquence immédiate du théorème précédent. En effet, si $P$ est un polynôme annulateur de $u$, alors $\ker P(u)=E$. $\blacksquare$

\subsubsection{Critère de diagonalisabilité}

 \begin{theo}[]

 $u\in\mathcal{L}(E)$ est diagonalisable si et seulement s'il existe un polynôme annulateur de $u$ scindé à racines simples sur $\mathbb{K}$.

\end{theo}

\paragraph{Preuve:}

\underline{Condition nécessaire:} 

$u$ est diagonalisable donc $E=\bigoplus\limits_{i=1}^r E_{\lambda_i}$ où les $ E_{\lambda_i}$ sont les sous-espaces propres de $u$.

Soit $P(X) = \prod\limits_{i=1}^r (X-\lambda_i)$, les $(X-\lambda_i)$ sont premiers entre eux deux à deux, donc $\ker P(u)= \bigoplus\limits_{i=1}^r \ker(u-\lambda_i \Id)$.

On obtient alors $\ker P(u)=E$ et donc $P(u)=0$.
\\
\\
\underline{Condition suffisante:} 

Supposons qu'il existe $P = \alpha (X-\lambda_1)\cdots (X-\lambda_p)$ scindé à racines simples qui annule $u$. Dans ce cas $\ker P(u) = \bigoplus\limits_{i=1}^p (u-\lambda_i \Id)$.

Or $\ker P(u)=E$, donc $E= \bigoplus\limits_{i=1}^p (u-\lambda_i \Id)$.

Soit alors $J=\{ 1\leq i \leq p ; \ker (u-\lambda_i) \neq \{0\}\}$. Si $i\in J$ alors $\lambda_i$ est valeur propre de $u$ et on a alors $E= \bigoplus\limits_{i\in J} (u-\lambda_i \Id)$ ce qui signifie que $u$ est diagonalisable. $\blacksquare$

\subsubsection{Restriction aux sous-espaces stables}

 \begin{prop}[]

Soit $u \in \mathcal{L}(E)$ diagonalisable et $F$ un sous-espace vectoriel de $E$ stable par $u$, alors l'endomorphisme induit par $u$ à $F$ est diagonalisable

\end{prop}

\paragraph{Preuve:}

On note $u_F$ l'endomorphisme induit par $u$ à $F$. Si $u$ est diagonalisable alors il existe $P\in\mathbb{K}[X]$ scindé à racines simples tel que $P(u)=0$. Or $F \subset E$ et $u(F) \subset F$, donc $P(u_{F})=0$. $u_{F}$ est alors diagonalisable. $\blacksquare$

\subsection{Théorème de Cayley-Hamilton}

\subsubsection{Sous-espace cyclique et matrice compagnon}

 \begin{prop}[]

Soit $P(X) = (-1)^n (X^n + a_{n-1} X^{n-1} + ...+ a_1 X + a_0)$ .
\\
Soit  $A=\left ( \begin{array}{cccc}
0 & ..... &0& -a_0      \\
 1&..... &.&-a_1                      \\
0 & .....& .  & .        \\
.&......&0&-a_{n-2} \\
0&....0&1&-a_{n-1}
\end{array}\right)\in\mathcal{M}_n(\mathbb{K})$ 

appelée matrice compagnon de $P$.
\\
Alors on a : $\chi_A(X)=P(X)$.

\end{prop}

\paragraph{Preuve:}

Par récurrence sur $n\geq 1$.

\underline{Initialisation:} 

Soit $A=(-a_0)$. Dans ce cas $\chi_A(X)= -a_0-X$ et $P(X)= (-1) (X+a_0) = -a_0-X$. Donc $\chi_A(X)=P(X)$ et la propriété est alors vraie pour $n=1$.

\underline{Hérédité:}

Supposons la propriété vraie jusqu'à l'ordre $n$ et démontrons qu'elle est alors vraie à l'ordre $n+1$.

Soit alors $A=\left ( \begin{array}{cccc}
0 & ..... &0& -a_0      \\
 1&..... &.&-a_1                      \\
0 & .....& .  & .        \\
.&......&0&-a_{n-1} \\
0&....0&1&-a_{n}
\end{array}\right)\in\mathcal{M}_{n+1}(\mathbb{K})$ 

\bigskip

On a $\chi_A(X)=\det(A -XI_{n+1})= \begin{vmatrix}
-X & ..... &0& -a_0      \\
 1&..... &.&-a_1                      \\
0 & .....& .  & .        \\
.&......&-X&-a_{n-1} \\
0&....0&1&-X-a_{n}
\end{vmatrix}$

En développant par rapport à la première ligne on obtient :

$\chi_A(X)= -X \begin{vmatrix}
-X & ..... &0& -a_1      \\
 1&..... &.&-a_2                      \\
0 & .....& .  & .        \\
.&......&-X&-a_{n-1} \\
0&....0&1&-X-a_{n}
\end{vmatrix} -(-1)^{n+2} a_0 \begin{vmatrix}
1 &-X & 0 &   .&.&0  \\
 0 &1&-X& 0&.&0                   \\
. &.& .  & . &.&.       \\
0& . & .& .  & 1&-X     \\
0 &.& . &. & .  &1 \\
\end{vmatrix}$

\bigskip

$\chi_A(X)=-X\Bigl((-1)^n(X^n+a_nX^{n-1}+\ldots+a_2X+a_1)\Bigr)+(-1)^{n+1}a_0$
\newline

$\chi_A(X)=(-1)^{n+1} (X^{n+1}+a_nX^n+\ldots+a_2X^2+a_1X+a_0)$
\newline

La propriété est alors vraie à l'ordre $n+1$. $\blacksquare$

\subsubsection{Théorème de Cayley-Hamilton}

 \begin{theo}[Cayley-Hamilton]

Pour tout $u\in\mathcal{L}(E)$ \ \  $\chi_u(u)=0$.

\end{theo}

\paragraph{Preuve:}

Soit $x\in E$, $(x,u(x), \ldots ,u^n(x))$ est une famille qui contient $n+1$ éléments. Or $\dim E = n$, donc cette famille est liée.

Soit $p$ le plus petit entier tel que $(x , u(x), \cdots  , u^p(x))$ soit liée.

Dans ce cas, $(x , u(x) , \ldots ,u^{p-1}(x))$ est libre et il existe $a_0,a_1, \ldots ,a_{p-1}\in\mathbb{K}$ tels que $u^p(x)+a_{p-1}u^{p-1}(x)+\ldots+a_0x=0$.

On complète la famille libre en une base de $E$ et on écrit la matrice de $M$ de $u$ dans cette base. On obtient alors :

$M=\left(\begin{array}{c|c}
A & C \\
\hline
0 & B \\

\end{array}
\right) \in\mathcal{M}_n(\mathbb{K})$ avec $A=\left ( \begin{array}{cccc}
0 & ..... &0& -a_0      \\
 1&..... &.&-a_1                      \\
0 & .....& .  & .        \\
.&......&0&-a_{p-2} \\
0&....0&1&-a_{p-1}
\end{array}\right)\in\mathcal{M}_p(\mathbb{K})$ .

\bigskip

On a $\chi_M(X)=\det (M-XI_n)=\det (A-XI_p) \det(B-XI_{n-p})=\chi_A(X) \chi_B(X)$.

 Or $\chi_A(X)=P(X)$ avec $P(X)=(-1)^p(X^p+\sum\limits_{k=p-1}^0 a_kX^k)$ et donc $\chi_A(A)(x)=0$.

Ceci étant vrai pour tout $x\in E$, donc $\chi_M(M)=0$. $\blacksquare$

\begin{corollaire}

Si  $A\in GL_n(\mathbb{K})$ alors $A^{-1}\in\mathbb{K}[A]$
\end{corollaire}

\paragraph{Preuve :}

 $\chi_A(X)= (-1)^n X^n+\sum\limits_{k=0}^{n-1} a_kX^k$. D'après le théorème de Cayley-Hamilton $\chi_A(A)=0$ , c'est-à-dire :

$(-1)^nA^n+\ldots+a_1A+a_0I_n=0 \iff (-1)^nA^n+\ldots+a_1A=-a_0I_n  $

En factorisant par $A^{-1}$ on obtient : $ A^{-1}\Bigl ( (-1)^nA^n+\ldots+a_1A \Bigr )= -a_0A^{-1}$.

Or $A$ est inversible et $\det A= a_0$, donc $a_0\neq 0$.

On a alors $A^{-1} = \displaystyle\frac{-1}{a_0} A^{-1}\left ( (-1)^nA^n+\ldots+a_1A \right )$ .

Donc  $A^{-1} = \displaystyle\frac{-1}{a_0} \left ( (-1)^nA^{n-1}+\ldots+a_1I_n \right )$  qui est un polynôme en $A$ .  $\blacksquare$

\subsection{Le lemme  chinois} 

Rappelons que $\mathbb{K}[X]$ est un anneau principal puisque $\mathbb{K}$ est un corps. 

De cette propriété tout idéal $I$ de  $\mathbb{K}[X]$ est principal et donc $I=(P)$ pour un certain $P\in \mathbb{K}[X]$ et dans ce cas  $\mathbb{K}[X]/I$ serait un anneau.

\subsubsection{Lemme chinois} 

 \begin{theo}[Chinois]

Soit $k\in\mathbb{N}$ et $k\geq2$. Soient $ P_1, \ldots , P_k$ des éléments de $\mathbb{K}[X]$ premiers entre eux deux à deux. Alors :

$$\mathbb{K}[X]/(P_1\cdots P_k) \simeq \mathbb{K}[X]/(P_1) \times \cdots \times\mathbb{K}[X]/(P_k)$$

\end{theo}

\paragraph{Preuve:}

Pour démontrer le théorème ci-dessus je propose deux preuves qui sont toutes les deux instructives car elles peuvent être exploitées dans d'autres démonstrations utilisant les mêmes propriétés intrinsèques.

\newpage

\underline{Première méthode:}

\begin{methode}
Soient deux anneaux $A$ et $B$. Soit $I$ un idéal  de $A$. 

Pour démontrer que $A/I \simeq B$, il suffit de construire un morphisme surjectif  $\varphi : A \to B$ tel que $I =\ker\varphi$. 

Dans ce cas $A/\ker\varphi \simeq \ $Im $\varphi$, ce qui est équivalent à  $A/I \simeq B$.

\end{methode}

Démontrons alors la propriété par récurrence sur $k\in\mathbb{N}$ et $k\geq 2$.
\\
\\
\underline{Initialisation:} 

Soient $P_1$ et $P_2$ deux éléments de  $\mathbb{K}[X]$ premiers entre eux. Il existe alors, d'après Bezout,  $U,V \in \mathbb{K}[X]$ tels que $UP_1+VP_2=1$.

\bigskip

Soit alors le morphisme d'anneaux  \[ \varphi : \left\{ 
                                                                          \begin{array}{ccc}\mathbb{K}[X] &\to & \mathbb{K}[X]/(P_1) \times  \mathbb{K}[X]/(P_2) \\
                                                                                                         P&\mapsto&(P \mod  P_1, P  \mod  P_2)
                                                                          \end{array} \right .\]

\bigskip

On a $\ker \varphi = \{P \in \mathbb{K}[X]$ tel que $P_1| P $ et $P_2 | P \}$. Or $P_1$ et $P_2$ sont premiers entre eux donc $\ker \varphi = \{P \in \mathbb{K}[X]$ tel que $P_1P_2| P \}=(P_1P_2)$.

On obtient donc $ \mathbb{K}[X]/ \ker \varphi \simeq$ Im $ \varphi$. Soit alors  $ \mathbb{K}[X]/ (P_1P_2) \simeq$ Im $ \varphi$

Montrons à présent que $\varphi $ est surjective.

 Pour cela, on considère $(Q , R)\in \mathbb{K}[X]/(P_1) \times \mathbb{K}[X]/(P_2)$ et on considère $P \in \mathbb{K}[X]$ avec $P=URP_1+VQP_2$.

Or $UP_1 \equiv 1 \ ( \mod\  P_2)$ et $VP_2 \equiv 1 \ (\mod \ P_1)$.

Donc  $P \equiv Q \ (\mod \ P_1)$ et $P \equiv R \ (\mod \ P_2)$.

 On a alors $\varphi(P) =(Q , R)$ et donc $\varphi$ est surjective avec 

Im$\varphi =  \mathbb{K}[X]/(P_1) \times  \mathbb{K}[X]/(P_2)$.

Finalement on a :  $ \mathbb{K}[X]/(P_1P_2) \simeq  \mathbb{K}[X]/(P_1) \times  \mathbb{K}[X]/(P_2)$.
\\
\\
\underline{Hérédité:}

Soient $k\in\mathbb{N}$ et $k\geq 2$ et soient $P_1, \ldots ,P_k,P_{k+1} \in \mathbb{K}[X]$ premiers entre eux deux à deux.

Supposons que $\mathbb{K}[X]/(P_1\cdots P_k) \simeq \mathbb{K}[X]/(P_1)  \times \cdots \times\mathbb{K}[X]/(P_k)$.

 Puisque $(P_1 \cdots P_k)$ et $P_{k+1}$ sont premiers entre eux, d'après l'initialisation on obtient $\mathbb{K}[X]/(P_1\cdots P_{k+1}) \simeq \mathbb{K}[X]/(P_1\cdots P_k)  \times \mathbb{K}[X]/(P_{k+1})$.

Avec l'hypothèse de récurrence on obtient : $$\mathbb{K}[X]/(P_1\cdots P_k P_{k+1}) \simeq \mathbb{K}[X]/(P_1)  \times\cdots \times\mathbb{K}[X]/(P_k)\times\mathbb{K}[X]/(P_{k+1})$$ $\blacksquare$

\newpage

\underline{Deuxième méthode :} 

\begin{methode}
Cette méthode nous permettra de résoudre des systèmes de congruences dans $\mathbb{K}[X]$ et sera reprise dans  la preuve de la décomposition de Dunford, d'où son intérêt.

Cette fois-ci, on passe par une méthode directe en démontrant que :

 \[  \varphi : \left\{ 
                            \begin{array}{ccc}
                              \mathbb{K}[X]/(P_1\cdots P_k) &\to & \mathbb{K}[X]/(P_1) \times \cdots \times  \mathbb{K}[X]/(P_k) \\
                              P  \mod (  P_1\cdots P_k )&\mapsto&(P  \mod  P_1 , \ldots , P \mod P_k)
                             \end{array} 
                              \right. \] 

est un isomorphisme.

\end{methode}

En effet, on définit $Q_i = \prod\limits_{j=1, j\neq i}^k P_j$ pour $ 1 \leq i \leq k$, et dans ce cas les $Q_i$ seront premiers entre eux dans l'ensemble.

 Par le théorème de Bezout, il existe $U_1, \ldots , U_k \in \mathbb{K}[X]$ tels que $U_1Q_1+\ldots+ U_kQ_k=1$. 

Pour tout $(R_1, \ldots ,R_k)\in  \mathbb{K}[X]/(P_1) \times \cdots \times  \mathbb{K}[X]/(P_k)$, onconsidère $$P = U_1Q_1R_1+\ldots+ U_kQ_kR_k$$ qui vérifie bien $P \equiv R_i \ (\mod \ P_i)$ pour $1\leq i \leq k$. 

Donc $\varphi(P)=(R_1, \ldots ,R_k)$ et $\varphi$ est alors surjective.

Puis $\ker \varphi= \{P \in \mathbb{K}[X]/(P_1\cdots P_k)$ tel que $P_i | P$ pour tout $1\leq i \leq k \}$.

Puisque les $P_i$ sont premiers entre eux deux à deux, alors :

$\ker \varphi= \{P \in \mathbb{K}[X]/(P_1\cdots P_k)$ tel que $P_1\cdots P_k | P \}=\{0\}$.

Donc $\varphi$ est injective et donc $\varphi$ est un isomorphisme. $\blacksquare$

\subsubsection{Système de congruences} 

 \begin{prop}[]

Soit $k\in\mathbb{N}$ et $k\geq2$. Soient $ P_1 \cdots P_k$ des éléments de $\mathbb{K}[X]$ premiers entre eux deux à deux. Alors le système : 

\[ (\mathcal{S}) : \left\{ \begin{array}{ccc}
P\equiv & Q_1& \mod  P_1     \\

P\equiv & Q_2& \mod P_2    \\

\cdots\equiv&\cdots&\cdots \\

P\equiv & Q_k& \mod  P_k     \\      
\end{array} \right. \]
\\
admet une solution unique $\mod ( P_1 \cdots P_k)$.

\end{prop}

\paragraph{Preuve:}

L'unicité de la solution provient de l'isomorphisme $\varphi$ construit dans la preuve précédente. La solution du système de congruences est obtenue  par la méthode qui a été utilisée pour démontrer la surjectivité de $\varphi$. $\blacksquare$

\subsection{Suite de Sturm}

\subsubsection{Théorème de Sturm}

\begin{theo}[de Sturm]

Soit $P\in\mathbb{R}[X]$.  On considère la suite $(S_i)_{1\leq i \leq p+1}$ définie par $S_0=P , S_1=P'$ et pour $i \geq 1$, $S_{i+1}$ est l'opposé du reste de la division euclidienne de $S_{i-1}$ par $S_{i}$.

 C'est-à-dire $\deg S_{i+1} < \deg S_i$ et il existe $Q_i$ tel que $S_{i-1}=Q_iS_i-S_{i+1}$. 

On suppose $S_{p+1}=0$ et donc $S_p=P \wedge P$'.
\\
\\

 Pour $x\in\mathbb{R}$, soit $V(x)$ le cardinal de l'ensemble

 $\{(i,j), 0\leq i<j \leq p, S_i(x)S_j(x) < 0 \text { et } \forall k\in [\![  i-1,j-1]\!], S_k(x)=0\}$ qui est le nombre de changement de signes de la suite $S_0(x), \ldots ,S_p(x)$.
\\
\\

Soient $a < b \in \mathbb{R}$, non racines de $P$. Alors le nombre de racines dans l'intervalle $[a,b]$ est égal à $V(a)-V(b)$.

\end{theo}

\paragraph{Preuve:}

 $S_{i+1}$ est l'opposé du reste de la division euclidienne de $S_{i-1}$ par $S_i$. Ces divisions euclidiennes successives sont possibles tant que $S_i$ n'est pas nul. 

 Donc la suite $S_0, \ldots ,S_p,S_{p+1}$ est au signe près la suite des restes obtenus en calculant le pgcd de $P$ et $P'$ par l'algorithme d'Euclide.

On alors pour $0\leq i \leq p$ , pgcd$(S_i,S_{i+1})=S_p$.

\bigskip

\textbf{\underline{Première étape :}}
\\
On se ramène à un polynôme n'ayant que des racines simples en posant pour  

$0\leq i \leq p$ ; $T_i=\displaystyle\frac{S_i}{S_p}$.

Alors, pour tout  $0\leq i \leq p$ , $T_i \wedge T_{i+1} = 1$ et dans ce cas $T_i$ et $T_{i+1}$ n'ont pas de racines communes et $T_p=1$.

Les racines de $T_0$ sont celles de $P$ et elles sont toutes simples pour $T_0$. En effet, si  $\alpha$ est racine de $P$ d'ordre $m$, alors $\alpha$ est racine de $P'$ d'ordre $m-1$ et dans ce cas $(X-\alpha)^{m-1}$ est en facteur dans $S_p$.

Si $S_p(x)\neq 0$, $V(x)$ est aussi le nombre de changement de signe dans la suite $T_0(x), \ldots , T_p(x)$, qu'on note $V_1(x)$.

Donc $V(a)-V(b) = V_1(a)-V_1(b)$.

\bigskip

\textbf{\underline{Deuxième étape :}}
\\
\\
Montrons que $ V_1(a)-V_1(b)$ est égal au nombre de racines de $P$ ( ou de $T_0$) sur $[a,b]$.

Pour y parvenir, on distingue trois cas possibles.
\bigskip

\underline{Premier cas :}
\\
Soit $I \subset [a,b]$ tel que :  $\forall x\in I, \forall i \in  $[\![$  0, p$]\!]$ , T_i(x) \neq 0$

Alors pour tout $x\in I$, $V_1(x)$ est constante car les $T_i$ sont des polynômes qui ne s'annulent pas sur $I$ et donc par continuité  sont de signe constant.
\bigskip

\underline{Deuxième cas :} 
\\
Soit $\alpha\in ]a,b[$ une racine de $P$. Alors $P^2(\alpha)=0$ et par continuité de $P$ , $P^2(x) \geq 0$ sur un voisinage de $\alpha$.

Donc la fonction $P^2$ admet un minimum en $\alpha$, et  elle change de variation lorsque $x$ traverse $\alpha$.

Donc il existe $\epsilon >0$ tel que $(P^2)' =2PP'$ est négative sur $[\alpha-\epsilon,\alpha]$ et puis positive sur $[\alpha,\alpha+\epsilon]$.

Soit alors $h\neq 0$ tel que $]\alpha-|h|,\alpha+|h|[ \subset ]\alpha-\epsilon,\alpha+\epsilon[ \subset ]a,b[$. 

Dans ce cas, le signe de $P(\alpha+h)P'(\alpha+h)$ est celui de $h$ car $(S_p(\alpha+h))^2 >0$ si $P(\alpha+h)\neq 0$.

Donc après la traversée de $\alpha$, il y a un changement de signe de moins entre les termes d'indice 0 et 1.
\bigskip

\underline{Troisième cas :}
\\
soit $\alpha$ une racine de $T_i$ avec $i\geq 1$. On a $i<p$ car $S_p(\alpha)\neq 0$ puisque $T_0(\alpha)\neq 0$.

De $S_{i-1}=Q_iS_i - S_{i+1}$, on déduit $T_{i-1}=Q_iT_i - T_{i-1}$ puis $T_{i-1}(\alpha) = -T_{i+1}(\alpha)$.

Or $\alpha$ n'est pas racine de $T_{i-1}$ ni de $T_{i+1}$, sinon on aurait $T_0(\alpha)=0$. Donc $T_{i-1}(\alpha)T_{i+1}(\alpha) < 0$ et par continuité dans un voisinage de $\alpha$, $T_{i-1}(x)T_{i+1}(x) < 0$.

Donc le nombre de changements de signe entre $T_{i-1}$ et $T_{i+1}$ ne change pas et il est toujours égal à 1 pour n'importe quel signe de $T_i(x)$ sur ce même voisinage.
\bigskip

Finalement les  trois cas possibles qu'on vient d'étudier constituent une partition de  l'intervalle $[a,b]$. Donc on  a bien $V(a)-V(b) = $ nombre de racines de $P$ dans l'intervalle $[a,b]$. $\blacksquare$

\subsubsection{Bornes de Cauchy}

\begin{prop}

Soit $P \in \mathbb{C}[X]$ avec $P = a_0 + a_1 X +\ldots+ a_k X^k $ tel que $a_k \neq 0$.

Alors, si $\alpha \in \mathbb{C}$ est une racine de $P$, on a :

$$|\alpha| \leq 1 + \max\limits_{1\leq i \leq k-1} \displaystyle\frac{|a_i|}{|a_k|}$$

\end{prop}

\goodbreak

\paragraph{Preuve :}

Soit $P\in\mathbb{C}[X]$

Soit $Q\in\mathbb{C}[X]$ tel que $Q(X) = \displaystyle\frac{P(X)}{a_k}$.

Soit $ bi = \displaystyle\frac{|a_i|}{|a_k|}$ et $B = \max\limits_{1\leq i \leq k-1}(b_i)$.
 
Le polynôme $Q$ est unitaire, a les mêmes racines que $P$ et vérifie :

Pour tout $x \in\mathbb{C}$ et $ |x| \neq 1$ : $|Q(x)| \geq |x|^k - B(|x|^{k-1} +\ldots+ 1) = |x|^k - B\displaystyle\frac{|x|^{k-1}-1}
{|x| - 1}$

 Si $|x| > 1 + B$, on a $1 > \displaystyle\frac{B}{|x|-1}$ et $|x|^k >\displaystyle\frac{|x|^kB}{|x|-1}$, 

d’où $|Q(x)| >\displaystyle\frac{B}{|x| - 1}(|x|^k - (|x|^k - 1))=
 \displaystyle\frac{B}{|x| -1} > 0$.

 Donc $x$ ne peut pas être racine de $Q$.

 Par suite si $\alpha$ est racine de $P$ alors $|\alpha|\leq 1+B$. $\blacksquare$

\paragraph{L'algorithme avec Xcas:}

Voici un programme qui prend trois arguments d'entrée : Le polynôme $P$ et les deux bornes d'un intervalle $[a,b]$.

Pour déterminer le nombre de racines (sans leur ordre de multiplicité) sur tout $\mathbb{R}$, il suffit de prendre $a= -M$ et $b=M$, où $M=1+B$ qui est la borne de Cauchy donnée dans la proposition précédente.

\parbox{12cm}{suite\_sturm (P,a,b):= \\
\{ local S0,S1,S2,T0,T1,T2,G,seqa,seqb,suba,subb,n,m,k,j,Vab; \\
  G:=gcd(P,diff(P,x));  \\
  S0:=P;  \\
  T0:=P/G;  \\
  seqa:=[subs(T0,x=a)];  \\
  seqb:=[subs(T0,x=b)];  \\
  S1:=diff(P,x);  \\
  T1:=S1/G;  \\
  suba:=subs(T1,x=a);  \\
  subb:=subs(T1,x=b);  \\
  if (suba{\tt\symbol{60}}{\tt\symbol{62}}0) seqa:=append(seqa,suba); ;  \\
  if (subb{\tt\symbol{60}}{\tt\symbol{62}}0) seqb:=append(seqb,subb); ;  \\
  while(S1{\tt\symbol{60}}{\tt\symbol{62}}0)\{ \\
      S2:=-rem(S0,S1);  \\
      T2:=S2/G;  \\
      suba:=subs(T2,x=a);  \\
      subb:=subs(T2,x=b);  \\
      if (suba{\tt\symbol{60}}{\tt\symbol{62}}0) seqa:=append(seqa,suba); ;  \\
      if (subb{\tt\symbol{60}}{\tt\symbol{62}}0) seqb:=append(seqb,subb); ;  \\
      S0:=S1;  \\
      S1:=S2;  \\
    \};; ;  \\
  n:=0;  \\
  m:=0;  \\
  for (k:=0;k{\tt\symbol{60}}(size(seqa)-1);k++) if ((seqa[k]*seqa[k+1]){\tt\symbol{60}}0) n:=n+1; ; ;  \\
  for (j:=0;j{\tt\symbol{60}}(size(seqb)-1);j++) if ((seqb[j]*seqb[j+1]){\tt\symbol{60}}0) m:=m+1; ; ;  \\
  Vab:=n-m;  \\
  return(Vab);  \\
\} }

\newpage

\section{Décomposition de Dunford} 

\subsection{Codiagonalisation} 

\begin{prop}[]

Soient $u$ et $v$ deux endomorphismes diagonalisables de $\mathcal{L}(E)$ qui commutent. Alors il existe une base commune de diagonalisation.

\end{prop}

\paragraph{Preuve:}

Supposons que  $u$ et $v$ sont deux endomorphismes diagonalisables de $\mathcal{L}(E)$ qui commutent.

Montrons d'abord que tout sous-espace propre de $u$ est stable par $v$.

En effet, soit $E_{\lambda}$ un sous-espace propre de $u$ et soit $x \in E_{\lambda}$.

$u(v(x)) = v(u(x)) = v(\lambda x) = \lambda v(x)$. Donc $v(x) \in  E_{\lambda}$.

Maintenant, considérons $\lambda_1, \ldots , \lambda_k$ les valeurs propres de $u$ et $ E_{\lambda_1}, \ldots , E_{\lambda_k}$ les sous-espaces propres associés.

Pour tout $1\leq i \leq k$, $ E_{\lambda_i}$ est stable par $v$ et la restriction de $v$ à $ E_{\lambda_i}$ induit un endomorphisme diagonalisable de $ E_{\lambda_i}$.

Il existe donc une base $B_i$ de $ E_{\lambda_i}$ de vecteurs propres de $v$ et aussi de $u$ puisque $u_{| E_{\lambda_i}} = \lambda_i \Id_{ E_{\lambda_i}}$.

On sait que $E=\bigoplus\limits_{i=1}^k  E_{\lambda_i}$ et donc $B=(B_1, \ldots ,B_k)$ est une base de diagonalisation commune de $u$ et $v$. $\blacksquare$

\begin{corollaire}

Soient $u$ et $v$ deux endomorphismes diagonalisables de $\mathcal{L}(E)$ qui commutent. Alors $u+v$ est diagonalisable.

\end{corollaire}

\paragraph{Preuve :}

Si  $u$ et $v$ deux endomorphismes diagonalisables de $\mathcal{L}(E)$ qui commutent alors il existe une base $B$ de veteurs propres communs à $u$ et $v$. La base $B$ est aussi une base de vecteurs propres de $u+v$. Donc $u+v$ est aussi digonalisable. $\blacksquare$

\subsection{Somme de nilpotents qui commutent} 

\begin{prop}[]

Soient $n$ et $n'$ deux endomorphismes nilpotents de  $\mathcal{L}(E)$ qui commutent. Alors $n+n'$ est nilpotent.

\end{prop}

\paragraph{Preuve:}

Soient $n$ et $n'$ deux endomorphisme nilpotents qui commutent. Soient alors $r$ et $r'$ les indices de nilpotence respectifs de $n$ et $n'$.

$(n+n')^{r+r'} = \displaystyle\sum\limits_{k=0}^{r+r'} \displaystyle\binom{r+r'}{k} n^k n'^{r+r'-k}$.

Si $k < r $ alors $ r+r'-k > r'$ et dans ce cas $n'^{r+r'-k} = 0$.

Si $k \geq r$ alors $n^k =0$.

Donc tous les termes de la somme sont nuls et  on aura  toujours : $(n+n')^{r+r'}=0$. $\blacksquare$

\begin{corollaire}

Soit $r\in\mathbb{N}$ et $r\geq 2$. Soient $n_1, \ldots , n_r$ des endomorphismes nilpotents qui commutent deux à deux. Alors $\sum\limits_{k=1}^r n_k$ est nilpotent.

\end{corollaire}

\paragraph{Preuve :}

Par récurrence sur  $r\in\mathbb{N}$ et $r\geq 2$.

Pour $r=2$, si $n_1$ et $n_2$ sont nilpotents qui commutent, alors $ n_1+n_2$ est nilpotent d'après la proposition précédente.

Supposons que la propriété reste vraie au rang $r$.

Soient $n_1, \ldots , n_r,n_{r+1}$ des endomorphismes nilpotents qui commutent deux à deux. D'après l'hypothèse de récurrence , $n=n_1+\ldots+n_r$ est nilpotent.

De plus $n$ et $ n_{r+1}$ sont aussi nilpotents qui commutent, d'après l'initialisation $n+n_{r+1}$  est nilpotent. $\blacksquare$

\subsection{Diagonalisable et nilpotent} 

\begin{prop}[]

Soit $u$ un endomorphismes diagonalisable de $\mathcal{L}(E)$. Si $u$ est nilpotent alors $u=0$.

\end{prop}

\paragraph{Preuve:}

Soit $u$ un endomorphisme diagonalisable et nilpotent. Soit $r$ l'indice de nilpotence de $u$.

Donc $\mu_u(X)=X^r$. Or $u$ est diagonalisable et donc $\mu_u$ est scindé à racines simples, c'est-à-dire que $r=1$. 
Dans ce cas $\mu_u(X)=X$. Puisque $\mu_u(u)=0$, alors $u=0$. $\blacksquare$

\subsection{Théorème de Dunford-Schwarz} 

\begin{theo}[Dunford-Schwarz]

Soit  $u$ un endomorphisme de $\mathcal{L}(E)$ dont le  polynôme caractéristique est scindé sur $\mathbb{K}$. Alors il existe un unique couple $(d,n) \in\mathcal{L}(E)^2$ vérifiant :

1) $u = d + n$ avec $d$ diagonalisable sur $\mathbb{K}$ et $n$ nilpotent.

2) $dn = nd$.

3) $d,n \in \mathbb{K}[u]$.

\end{theo}

\paragraph{Preuve:}

\underline{Première étape:} 

Supposons que $(-1)^n \chi_u(X)=\prod\limits_{i=1}^r (X-\lambda_i)^{\alpha_i}$ avec $\lambda_i\neq \lambda_j$ pour $i\neq j$ et $\alpha_i \geq 1$.

Par le théorème de Cayley-Hamilton $\chi_u(u)=0$ et donc d'après le théorème de décomposition des noyaux : $\ker \chi_u(u) = E = \bigoplus\limits_{i=1}^r \ker (u-\lambda_i \Id)^{\alpha_i}$.

On note $N_i=\ker (u-\lambda_i \Id)^{\alpha_i}$ et $P_i=(X-\lambda_i)^{\alpha_i}$.

Soit $Q_i=\prod\limits_{k=1,k\neq i}^rP_k$, les $Q_i$ sont alors premiers entre eux dans leur ensemble et d'après le théorème de Bezout, il existe $U_1, \ldots ,U_r\in\mathbb{K}[X]$ tels que $U_1Q_1+\ldots+U_rQ_r=1$.

On a alors $U_1Q_1(u)+\ldots+U_rQ_r(u)=\Id$ et pour tout $x\in E$ on a 

$x=\sum\limits_{i=1}^r U_i(u) \circ Q_i(u)(x)$.

Notons $\pi_i= U_i(u) \circ Q_i(u)$. Pour $i\neq j$, $\pi_i \circ \pi_j= (U_iU_j)(u) \circ (Q_iQj)(u)$.

Or $\chi_u$ divise $Q_iQ_j$ et $\chi_u(u)=0$. Donc pour $i\neq j$, $\pi_i \circ \pi_j = 0$.

D'autre part $\pi_i=\pi_i \circ \Id = \pi_i \circ \sum\limits_{j=1}^r \pi_j = \sum\limits_{j=1}^r (\pi_i \circ \pi_j)$.

Donc $\pi_i=\pi_i \circ \pi_i$ et dans ce cas \underline{$\pi_i$ est un projecteur de $\mathbb{K}[u]$}.
\bigskip

Soit $x\in$ Im$\pi_i$. Il existe alors $y\in E$ tel que $x=\pi_i(y)=U_i(u)\circ Q_i(u)(y)$.

Dans ce cas $P_i(u)(x)=U_i(u)\circ (P_iQ_i)(u)(y)=U_i(u)\circ \chi_u(u)(y)=0$ (car $P_iQ_i=\chi_u$ et $\chi_u(u)=0$).

On a alors Im$\pi_i \subset N_i$.

Soit $x\in N_i$. On sait que $x=\sum\limits_{j=1}^r U_j(u)\circ Q_j(u)(x)$ et donc $x=U_i(u)\circ Q_i(u)(x)=\pi_i(x)$.

Alors $x\in\text{Im} \pi_i$ et on a dans ce cas \underline{Im$\pi_i = N_i$}.

Soit $x\in N_j$ pour $j\neq i$, $\pi_i(x)=U_i(u)\circ Q_i(u)(x)=0$. Donc $\bigoplus\limits_{j=1,j\neq i}^r N_j \subset \ker \pi_i$.

Puis si $x\in \ker \pi_i$, $x=\sum\limits_{j=1}^r \pi_j(x)=\sum\limits_{j=1,j\neq i}^r \pi_j(x)$.

Donc $x\in \bigoplus\limits_{j=1,j\neq i}^r N_j$ et on a alors \underline{$\ker \pi_i =\bigoplus\limits_{j=1,j\neq i}^r N_j$}.

\begin{retenir}

Pour $1\leq i \leq r$, $\pi_i$ est appelé projecteur spectral associé à la valeur propre $\lambda_i$ et vérifie les propriétés suivantes:

\begin{itemize}
\item {$\pi_i $ est un projecteur sur $N_i$ parallèlement à $\bigoplus\limits_{j=1,j\neq i}^r N_j$}.
\item {$\pi_i$ est un polynôme en $u$, c'est-à-dire que $\pi_i \in \mathbb{K}[u]$}.
\end{itemize}

\end{retenir}

\underline{Deuxième étape:}

On pose $d=\sum\limits_{i=1}^r \lambda_i \pi_i \in \mathbb{K}[u]$.

Or $E=\bigoplus\limits_{i=1}^r \ker (u-\lambda_i \Id)^{\alpha_i} = \bigoplus\limits_{i=1}^r N_i$.

Montrons alors que dans ce cas $N_i$ est un sous-espace propre de $d$.

En effet pour $x\in N_i$, $d(x)=\sum\limits_{j=1}^r \lambda_j\pi_j(x)=\lambda_i\pi_i(x)=\lambda_i x$.

Donc $E$ est somme directe des sous-espaces propres de $d$. $d$ est alors diagonalisable.

Posons $n=u-d$. On a dans ce cas $n=u \circ \Id - d = u \circ \left ( \sum\limits_{i=1}^r \pi_i \right )-\sum\limits_{i=1}^r \lambda_i \pi_i$.

On a alors $n = \sum\limits_{i=1}^r (u-\lambda_i \Id)\circ \pi_i$.

Soit $q =\max\limits_{1\leq i \leq r} \alpha_i$. On a alors :

$n^q=\sum\limits_{i=1}^r (u-\lambda_i \Id)^q \circ \pi_i=0$, car $\chi_u$  divise $(X-\lambda_i)^q U_iQ_i$.

Donc $n$ est un élément de $\mathbb{K}[u]$ qui est nilpotent. Puisque $d\in\mathbb{K}[u]$, on a $dn=nd$.

\begin{retenir}

\begin{itemize}
\item  {$d=\sum\limits_{i=1}^r \lambda_i \pi_i \in \mathbb{K}[u]$ et $d$  est diagonalisable}.
\item {$n=u-d \in\mathbb{K}[u]$ et $n$ est nilpotent}.
\item {$d$ et $n$ commutent.}
\end{itemize}

\end{retenir}

\underline{Troisième étape :}

Montrons l'unicité de la décomposition $u =d +n$.

Supposons alors qu'il existe un couple $(d',n') \in\mathbb{K}[u]^2$ vérifiant $u=d'+n'$.

On aura dans ce cas $d+n=d'+n'$, c'est-à-dire $d - d' = n' - n$.

Or $d$ et $d'$ sont diagonalisables qui commutent et donc ils sont diagonalisables dans une base commune $B$. $B$ est alors une base de diagonalisation de $d-d'$.

Puis $n$ et $n'$ sont nilpotents qui commutent et donc $ n' - n$ est aussi nilpotent. 

Or un endomorphisme qui est à la fois diagonalisable et nilpotent est obligatoirement nul, donc $d -d' = n' - n =0$. Finalement  $d=d'$ et $n=n'$. $\blacksquare$

\begin{corollaire}

Si $E$ est $\mathbb{C}$-espace vectoriel de dimension $n\in\mathbb{N}^*$, alors tout endomorphisme de $\mathcal{L}(E)$ admet une décomposition de Dunford.

\end{corollaire}

\paragraph{Preuve :}

En effet, $\mathbb{C}[X]$ est algébriquement clos, donc tout polynôme de $\mathbb{C}[X]$ est scindé.
 En particulier, tout polynôme annulateur de $u\in\mathcal{L}(E)$ est scindé.

D'où le résultat d'après le théorème de Dunford-Schwarz. $\blacksquare$

\begin{remarque}

Le théorème de Dunford-Schwarz reste valable dans le cas où $u$ admet un polynôme annulateur scindé sur $\mathbb{K}$.
 En particulier on peut énoncer le théorème avec le polynôme minimal de $u$, puisque $\mu_u$ divise $\chi_u$ et $\mu_u$ et $\chi_u$ ont les mêmes racines dans $\mathbb{K}$.

\end{remarque}

\subsection{Exemple de décomposition} 

La décomposition de Dunford par le calcul des projecteurs spectraux n'est réalisable dans la pratique que lorsqu'on connaît les valeurs propres. Pour y parvenir, on décompose en éléments simples $\displaystyle\frac{1}{P}$ où $P$ est un polynôme annulateur de $u$ .

Soit  $P=(-1)^n \chi_u=\prod\limits_{i=1}^r(X-\lambda_i)^{\alpha_i}$ qui est annulateur de $u$ .

 $\displaystyle\frac{1}{\chi_u} = \sum\limits_{i=1}^r \left(\sum\limits_{k=1}^{\alpha_i} \displaystyle\frac{a_{i,k}}{(X-\lambda_i)^k}\right)=\sum\limits_{i=1}^r \left(\sum\limits_{k=1}^{\alpha_i} \displaystyle\frac{a_{i,k}(X-\lambda_i)^{\alpha_i-k}}{(X-\lambda_i)^{\alpha_i}}\right)$.

On pose $U_i=\sum\limits_{k=1}^{\alpha_i} a_{i,k}(X-\lambda_i)^{\alpha_i-k}$, et on obtient  $\displaystyle\frac{1}{\chi_u}=\sum\limits_{i=1}^r\displaystyle\frac{U_i}{(X-\lambda_i)^{\alpha_i}}$.

Puis en multilpliant par $\chi_u$, on obtient : $1 =  \sum\limits_{i=1}^rU_iQ_i=\sum\limits_{i=1}^r P_i$ , où $Q_i=\prod\limits_{i\neq j}(X-\lambda_j)^{\alpha_j}$.

Dans ce cas chaque projecteur spectral est $\pi_i=P_i(u)$.

\bigskip

Pour cela, prenons un exemple. Soit $A=\begin{pmatrix} 
0&8&-4&8&-16&10\\
-8&14&-2&8&-18&8\\
-6&11&0&11&-25&13\\
2&7&-8&15&-21&17\\
2&1&-2&3&-3&5\\
0&0&0&0&0&4
\end{pmatrix}$

$$\chi_A(x)=x^6-30 x^5+368 x^4-2368 x^3+8448 x^2-15872 x+12288$$

$$\displaystyle\frac{1}{\chi_A}=\displaystyle\frac{1}{512(x-8)}+\displaystyle\frac{1}{-32(x-6)}+\displaystyle\frac{1}{8(x-4)^4}+\displaystyle\frac{3}{32(x-4)^3}+\displaystyle\frac{7}{128(x-4)^2}+\displaystyle\frac{15}{512(x-4)}$$.
\\
$P_1=\displaystyle\frac{\chi_A}{512(x-8)}=\displaystyle\frac{x^5-22x^4+192x^3-832x^2+1792x-1536}{512}$
\\
\\
$P_2=\displaystyle\frac{\chi_A}{-32(x-6)}=\displaystyle\frac{-x^5+24x^4-224x^3+1024x^2-2304x+2048}{32}$
\\
\\
$P_3=\chi_A\Bigl(\displaystyle\frac{3}{32(x-4)^3}+\displaystyle\frac{7}{128(x-4)^2}+\displaystyle\frac{15}{512(x-4)}\Bigr)$
\\
\\
$P_3=\displaystyle\frac{7x^4-142x^3+1032x^2-3232x+3840}{128}$
\\
\\
Enfin $D=8P_1(A)+6P_2(A)+4P_3(A)$ , ce qui donne :
\\
\\
 $D=\begin{pmatrix}
-12&20&-10&14&-16&16\\
-20&26&-8&14&-18&14\\
-26&31&-10&21&-25&23\\
-18&27&-18&25&-21&27\\
-6&9&-6&7&-3&9\\
0&0&0&0&0&4
\end{pmatrix}$ et $N=\begin{pmatrix}
12 & -12 & 6 & -6 & 0 & -6 \\
12 & -12 & 6 & -6 & 0 & -6 \\
20 & -20 & 10 & -10 & 0 & -10 \\
20 & -20 & 10 & -10 & 0 & -10 \\
8 & -8 & 4 & -4 & 0 & -4 \\
0 & 0 & 0 & 0 & 0 & 0
\end{pmatrix}$

\subsection{Une utilisation de la densité...}

Maintenant qu'on sait que toute matrice $M$ de $\mathcal{M}_n(\mathbb{C})$ admet une décomposition unique $D+N$ avec $D$ diagonalisable et $N$ nilpotente, on peut définir l'application

 \[ \varphi :\left \{
                   \begin{array}{ccc}
                    \mathcal{M}_n(\mathbb{C})&\to&\mathcal{D}_n(\mathbb{C})\\
                    M&\mapsto&D
                   \end{array} \right. \]

où $\mathcal{D}_n(\mathbb{C}) $ est l'ensemble des matrices diagonalisables sur~ $\mathbb{C}$.

$\varphi$ est-elle continue ? Pour répondre à cette question on rappelle tout d'abord quelques propriétés :

\begin{itemize}
\item{Toutes les normes sur $\mathcal{M}_n(\mathbb{C})$ sont équivalentes}

\item{Toute matrice de $\mathcal{M}_n(\mathbb{C})$ est trigonalisable}

\item{$\mathcal{D}_n(\mathbb{C}) $ est dense dans $\mathcal{M}_n(\mathbb{C}) $}
\end{itemize}

\begin{prop}

Pour $n\geq 2$, $\varphi$ n'est pas continue.

\end{prop}

\paragraph{Preuve :}

Soit $D\in\mathcal{D}_n(\mathbb{C}) $, alors $\varphi(D)=D$. Donc $\varphi$ est l'Identité sur $\mathcal{D}_n(\mathbb{C}) $.

Si $\varphi$ était continue, par densité $\varphi$ serait aussi l'Identité sur $\mathcal{M}_n(\mathbb{C}) $.

Or une matrice $N$  nilpotente non nulle n'est pas diagonalisable  et $\varphi(N)=N$ qui serait diagonalisable. Ceci est impossible, donc $\varphi$ n'est pas continue. $\blacksquare$

\section{Méthodes effectives} 

\subsection{Méthode de Newton-Raphson} 

On se place dans le cas où $\mathbb{K}=\mathbb{C}$ et onconsidère $E$ un $\mathbb{K}$-espace vectoriel de dimension $n\in\mathbb{N}^*$.

 Donc pour tout $u\in\mathcal{L}(E)$ , son polynôme caractéristique  $\chi_u$ est scindé, d'après le théorème de Dunford $u=d+n$ où $d$ est diagonalisable sur $\mathbb{C}$.

On a vu que lorsque le spectre de $u$ est connu, une décomposition de $\displaystyle\frac{1}{\chi_u}$ en éléments simples permet de trouver $d$ et $n$.

Cependant, lorsqu'on n'est pas en mesure de trouver les valeurs exactes des racines de $\chi_u$, la situation devient plus difficile puisque l'écriture de $d=\sum\limits_{i} \lambda_i\pi_i$ n'est plus possible.
\\
\\
L'Idée est de proposer un bon candidat $d$ qui sera annulé par un polynôme dans $\mathbb{C}[X]$ scindé à racines simples  puisque $d$ est diagonalisable.

Soit alors $\chi_u(X)=(-1)^n \prod\limits_{i=1}^r(X-\lambda_i)^{\alpha_i}$, par l'algorithme d'Euclide  on est capable de calculer le polynôme scindé à racines simples $P(X)=\prod\limits_{i=1}^r(X-\lambda_i)$.

En effet, $P=\displaystyle\frac{\chi_u}{\chi_u\wedge\chi'_u}$, et dans ce cas on cherche une solution $d$ vérifiant l'équation $P(d)=0$. 

L'Idée est de reprendre la méthode de résolution des équations numériques avec la méthode de Newton et de l'adapter à notre situation. On obtient alors une suite stationnaire vers un endomorphisme diagonalisable $d$.

Pour vérifier que $d$ correspond à la composante diagonalisable de Dunford il suffit de vérifier que $u-d$ est nilpotent et par le théorème de Dunford l'unicité est garantie.
\\
\\
Finalement on arrive à obtenir une méthode effective pour effectuer la décomposition de Dunford en se passant des valeurs propres. 
De plus cette méthode effective nous fournit un algorithme programmable permettant d'obtenir la décomposition souhaitée.

 \begin{lemme}[]

Si $U,N\in\mathcal{M}_n(\mathbb{K})$ tels que $U$ soit inversible et $N$ nilpotente qui commutent, alors $U-N$ est inversible.

\end{lemme}

\paragraph{Preuve :}

Supposons que $N$ est nilpotente d'indice $q\in\mathbb{N}^*$.

On a $(I_n - N) \sum\limits_{k=0}^{q-1} N^k = \sum\limits_{k=0}^{q-1} N^k - \sum\limits_{k=0}^{q-1} N^{k+1} = N^0 - N^q = I_n$.

Or $U - N = (I_n - NU^{-1}) U$ et $NU = UN$.

 On a alors $NU^{-1} = U^{-1}N$ et dans ce cas $(NU^{-1})^q = N^qU^{-q} = \huge 0_n$.

On en déduit que $\left ( \sum\limits_{k=0}^{q-1} (NU^{-1})^k \right ) (I_n - NU^{-1}) U = U$.

Donc $\left ( \sum\limits_{k=0}^{q-1} N^k U^{-k-1} \right ) (U - N) =I_n$.

Finalement $U - N$ est inversible. $\blacksquare$

\begin{lemme}[]

Soit $A \in \mathcal{M}_n(\mathbb{K})$ une matrice de polynôme caractéristique $\chi_A = \prod\limits_{i=1}^p (X-\lambda_i)^{\alpha_i}$. 

On pose $P =  \prod\limits_{i=1}^p (X-\lambda_i)$. Alors on a :

1) $P=\displaystyle\frac{\chi_A}{\chi_A \wedge \chi'_A}$ .

2)  $P'(A)$ est inversible.

3)  $P'(A)^{-1}$ est un polynôme en $A$.

\end{lemme}

\paragraph{Preuve:}

1) Les racines de $\chi_A$ sont les $\lambda_i$ de multiplicité $\alpha_i$, donc les racines de $\chi_A'$ sont les $\lambda_i$ de multiplicité $\alpha_i-1$.

Donc $\chi_A \wedge \chi_A' = \prod\limits_{i=1}^p (X - \lambda_i)^{\alpha_i-1}$.

On a alors $\displaystyle\frac{\chi_A}{\chi_A \wedge \chi_A'} =  \prod\limits_{i=1}^p (X - \lambda_i) = P(X)$. $\blacksquare$

2) $P$ et $P'$ sont premiers entre eux car ils n'ont pas de racines communes. D'après le théorème de Bezout il existe $U,V\in\mathbb{K}[X]$ tels que $UP+VP'=1$.

Donc $U(A)P(A)+V(A)P'(A)=I_n$ ce qui donne  $V(A)P'(A)=I_n - U(A)P(A) $.

Or d'après le théorème de Cayley-Hamilton $\chi_A(A)=0$ et donc pour $\alpha= \max_{1\leq i \leq p} \alpha_i$ on a $P(A)^{\alpha}=0$.

D'autre part $U(A)$ et $P(A)$ sont des polynômes en $A$, donc commutent et on a dans ce cas $(U(A)P(A))^{\alpha}=0$, c'est-à-dire que $U(A)P(A)$ est nilpotent.

On a alors  $I_n$  inversible et $U(A)P(A)$  nilpotent qui commutent, et donc d'après le lemme 1 , $P'(A)$ est inversible. $\blacksquare$

3) $P'(A)$ est inversible, d'après le corollaire du théorème de Cayley-Hamilton :

  $P'(A)^{-1}\in\mathbb{K}[A]$. $\blacksquare$

\begin{lemme}[]

Pour tout $Q \in \mathbb{K}[X]$ , il existe $R \in \mathbb{K}[X,Y]$ tel que :

$Q(X+Y)=Q(X)+YQ'(X)+Y^2R(X,Y)$

\end{lemme}

\paragraph{Preuve:}

Il suffit de démontrer la propriété pour un monôme. Soit alors $Q(X)= X^k$.

$Q(X+Y) = (X+Y)^k = \sum\limits_{i=0}^k \displaystyle\binom{k}{i} X^{k-i} Y^i = X^k + Y X^{k-1} + \sum\limits_{i=2}^k \displaystyle\binom{k}{i} X^{k-i} Y^i$.

Donc $Q(X+Y) = X^k + Y (kX^{k-1}) + Y^2 \left (\sum\limits_{i=2}^k \displaystyle\binom{k}{i} Y^{i-2} X^{k-i} \right )$

On a bien $Q(X+Y) = Q(X)+YQ'(X)+Y^2R(X,Y)$. $\blacksquare$

\begin{theo}[ Méthode de Newton-Raphson] 

Soit $A   \in\mathcal{M}_n(\mathbb{K})$ avec $n\in\mathbb{N}^*$, tel que $\chi_A = \prod\limits_{i=1}^p (X-\lambda_i)^{\alpha_i}$ et $P(X)= \prod\limits_{i=1}^p (X-\lambda_i)$  avec pour tout $1\leq i,j \leq r$, $\lambda_i \in \mathbb{C}$  , $\alpha_i\geq1$ et $\lambda_i \neq \lambda_j$ pour $i\neq j$. 
On définit la suite  \[ A_m= \left\{
                                             \begin{array}{cc}
                                             A_0=A &     \\
                          
                                              A_{m+1}=A_m-P(A_m)(P'(A_m))^{-1}  &      
                                            \end{array} \right. \]
\\
Alors la suite $A_m$ est bien définie et elle est stationnaire à partir d'un certain rang vers une matrice diagonalisable D qui réalise la décomposition de Dunford $A=D+N$.

\end{theo}

\paragraph{Preuve:}

\underline{\textbf{Première étape :}}

\bigskip

 Montrons par récurrence que la suite est bien définie et qu’elle est constituée de polynômes
en $A$. Autrement dit, pour tout $n\in\mathbb{N}$, $A_n$ est un polynôme en $A$ ,
 $P'(A_n)$ est inversible et $P'(A_n)^{-1}$ est aussi polynôme en $A$.

\bigskip

\underline{Initialisation :}  $A_0 = A$ est bien polynôme en $A$. Puis d'après le lemme 2,  $P'(A)$ est inversible et $P'(A)^{-1}$ est aussi polynôme en $A$.

\bigskip

\underline{Hérédité :}  Supposons pour tout $k \leq n$, $A_k$ est un polynôme en $A$
et $P'(A_k)$ est inversible. Alors $A_{n+1}$ est un polynôme en $A$ par hypothèse de récurrence et
parce que l’inverse d’une matrice est un polynôme en la matrice par le théorème de Cayley-Hamilton(même argument utlisé dans la preuve du 3) du lemme 2). 

Reste seulement à prouver que $P'(A_{n+1})$ est inversible pour conclure la récurrence.

Par la formule de Taylor $P'(A_{n+1})-P'(A_n) = (A_{n+1}-A_n)Q(A_n) $ pour un certain polynôme $Q$ soit $P'(A_{n+1}) = P'(A_n) - P(A_n)P'(A_n)^{-1}Q(A_n)$.

Or $P'(A_n)$ est un polynôme en $A$ qui est inversible par hypothèse de récurrence  .

Donc $ P(A_n)P'(A_n)^{-1}Q(A_n)$ est un polynôme en $A_n$ (qui est polynôme en $ A $ par hypothèse de récurrence). 
Alors $P'(A_n)$ et $P(A_n)P'(A_n)^{-1}Q(A_n)$  commutent. 
\\
\\
\underline{\textbf{Deuxième étape :}} 

\bigskip

Montrons par une autre récurrence que $P(A_n) \in P(A)^{2^n} \mathbb{K}[A]$.
\\
\\
\underline{Initialisation :}  On a $2^0=1$ et $P(A_0)=P(A)=P(A)^{2^0}$
\\
\\
\underline{Hérédité :}  En appliquant le lemme 3, il vient : 
$P(A_{n+1}) = P(A_n + Y ) = Y^2Q(A_n, Y )$ pour un certain polynôme $Q$, car $ P(A_n) + Y P'(A_n) = 0$ par définition de $Y$ . 

Ainsi $P(A_{n+1}) = P(A_n)^2P'(A_n)^{-2}Q(An, Y ) = (P(A)^{2^n})^2S(A)$ pour un certain polynôme $S$ par hypothèse de récurrence, 
soit $P(A_{n+1})\in P(A)^{2^{n+1}}\mathbb{K}[A]$ ce qui conclut la récurrence.

Revenons à la première récurrence :

$ P(A_n)P'(A_n)^{-1}Q(A_n)$ est bien nilpotent puisque $P(A)^r =0$ pour tout $r$ tel que $r=2^n > \max(\alpha_i)$. 

Le lemme 1 montre donc que $P'(A_{n+1})$ est inversible.

\bigskip

\underline{\textbf{Troisième étape :}}

\bigskip

 Montrons que la suite converge bien vers  $D$, la composante digonalisable de Dunford. 

En effet on vient de montrer que la suite est stationnaire vers $A_m$ dès que $2^m > \max(\alpha_i)$. 
Ainsi la suite converge en un nombre fini d’étapes, et il en faut seulement $E( \displaystyle\frac{\log(\max(\alpha_i))}{\log(2)} )+1$. 

De plus la limite $A_m = D$ vérifie $P(D)P'(D)^{-1} = 0$. 
 Or $P'(D)^{-1}$ est inversible donc $P(D) = 0$, et dans ce cas $D$ est diagonalisable car $P$ est scindé à racines simples. 

\bigskip

\underline{\textbf{Quatrième étape :}}

\bigskip

$N := A - D = A_0 - A_n = (A_0 - A_1) +\ldots+ (A_{n-1} - A_n)$ est une somme de nilpotents qui commutent 
 (comme on l’a vu lors de l’hérédité pour la première récurrence) donc $N$ est nilpotente. 

De plus $N$ est un polynôme en $A$ car $D$ l’est, donc $N$ commute avec $D$. 

 L’unicité dans la décomposition de Dunford nous permet alors de conclure. $\blacksquare$

\newpage

\paragraph{L'alogorithme en Xcas}   

On applique la méthode ci-dessus avec un programme qui prend la matrice $A$ à décomposer pour argument d'entrée et nous fournit à la sortie  $D$ et $N$.

\bigskip

\parbox{12cm}{dunford\_effective (A):= \\
\{ local chi,P,Q,An,Ann; \\
  chi:=pcar(A,x);  \\
  P:=simplify(chi/gcd(chi,diff(chi,x)));  \\
  Q:=diff(P,x);  \\
  An:=A;  \\
  Ann:=An-(subs(P,x=An))/(subs(Q,x=An));  \\
  while(An{\tt\symbol{60}}{\tt\symbol{62}}Ann)\{ \\
      An:=Ann;  \\
      Ann:=An-(subs(P,x=An))/(subs(Q,x=An));  \\
    \};; ;  \\
  return(An,A-An);  \\
\} }

\bigskip

\paragraph{Exemple :}

Onconsidère la matrice 

$A=\left(\begin{array}{cccccccc}
1 & 163 & 1 & 0 & 162 & 165 & 163 & -1 \\
0 & 209 & -2 & 0 & 210 & 209 & 208 & 2 \\
0 & -86 & -1 & 1 & -87 & -87 & -87 & 1 \\
1 & -22 & 4 & -1 & -23 & -20 & -20 & -4 \\
1 & 155 & 4 & 0 & 153 & 157 & 155 & -3 \\
0 & -111 & -1 & 0 & -110 & -112 & -110 & 1 \\
-1 & -245 & -1 & 0 & -244 & -246 & -245 & 1 \\
-1 & -160 & -3 & 1 & -158 & -162 & -161 & 4
\end{array}\right) \in\mathcal{M}_8(\mathbb{R})$
\\
\\
$$\chi_A(x) = x^{8}-8 x^{7}+6 x^{6}+16 x^{5}+97 x^{4}+216 x^{3}+264 x^{2}+288 x+144$$

On peut voir le polynôme caractéristique  $\chi_A$ comme étant un polynôme de $\mathbb{C}[X]$ qui est scindé sur $\mathbb{C}$.

 Donc la méthode de Newton-Raphson s'applique sur la matrice $A$ pour nous donner la décomposition de Dunford $A=D+N$ où $D$ est diagonalisable sur $\mathbb{C}$ .

On peut alors se poser la question suivante : $D$ est-elle diagonalisable sur $\mathbb{R}$ ?

Je propose une méthode de vérification par les suites de Sturm dans la dernière partie pour répondre à cette question.
Revenons à notre exemple, avec Xcas on exécute dunford\_effective(A) et on obtient :

\bigskip

$D=\left(\begin{array}{cccccccc}
\frac{3719}{16686} & \frac{315275}{3708} & -\frac{117607}{33372} & -\frac{81239}{33372} & \frac{1219219}{16686} & \frac{471985}{5562} & \frac{2927425}{33372} & -\frac{227905}{16686} \\
&&&&&&&\\
-\frac{23185}{16686} & \frac{240737}{3708} & -\frac{357343}{33372} & -\frac{239051}{33372} & \frac{379130}{8343} & \frac{360811}{5562} & \frac{2330431}{33372} & -\frac{208553}{8343} \\
&&&&&&&\\
\frac{6833}{11124} & -\frac{75491}{1648} & \frac{8389}{44496} & \frac{70883}{44496} & -\frac{893137}{22248} & -\frac{338575}{7416} & -\frac{2108287}{44496} & \frac{91001}{11124} \\
&&&&&&&\\
\frac{43421}{33372} & \frac{421813}{14832} & \frac{1135663}{133488} & \frac{451145}{133488} & \frac{2223845}{66744} & \frac{638927}{22248} & \frac{3620027}{133488} & \frac{200327}{33372} \\
&&&&&&&\\
\frac{7133}{11124} & \frac{443321}{4944} & \frac{59371}{44496} & -\frac{45403}{44496} & \frac{1694753}{22248} & \frac{659375}{7416} & \frac{4070207}{44496} & -\frac{166171}{11124} \\
&&&&&&&\\
\frac{1022}{2781} & -\frac{103061}{2472} & \frac{64301}{22248} & \frac{65155}{22248} & -\frac{329951}{11124} & -\frac{155297}{3708} & -\frac{968783}{22248} & \frac{40370}{2781} \\
&&&&&&&\\
\frac{653}{1854} & -\frac{266801}{2472} & \frac{47323}{7416} & \frac{38129}{7416} & -\frac{320869}{3708} & -\frac{132163}{1236} & -\frac{835825}{7416} & \frac{47741}{1854} \\
&&&&&&&\\
-\frac{5765}{11124} & -\frac{428887}{4944} & -\frac{7621}{44496} & \frac{101077}{44496} & -\frac{1597535}{22248} & -\frac{638045}{7416} & -\frac{3990233}{44496} & \frac{196861}{11124}
\end{array}\right) $

\bigskip

\[ N=\left(\begin{array}{cccccccc}
\frac{12967}{16686} & \frac{289129}{3708} & \frac{150979}{33372} & \frac{81239}{33372} & \frac{1483913}{16686} & \frac{445745}{5562} & \frac{2512211}{33372} & \frac{211219}{16686} \\
&&&&&&&\\
\frac{23185}{16686} & \frac{534235}{3708} & \frac{290599}{33372} & \frac{239051}{33372} & \frac{1372900}{8343} & \frac{801647}{5562} & \frac{4610945}{33372} & \frac{225239}{8343} \\
&&&&&&&\\
-\frac{6833}{11124} & -\frac{66237}{1648} & -\frac{52885}{44496} & -\frac{26387}{44496} & -\frac{1042439}{22248} & -\frac{306617}{7416} & -\frac{1762865}{44496} & -\frac{79877}{11124} \\
&&&&&&&\\
-\frac{10049}{33372} & -\frac{748117}{14832} & -\frac{601711}{133488} & -\frac{584633}{133488} & -\frac{3758957}{66744} & -\frac{1083887}{22248} & -\frac{6289787}{133488} & -\frac{333815}{33372} \\
&&&&&&&\\
\frac{3991}{11124} & \frac{322999}{4944} & \frac{118613}{44496} & \frac{45403}{44496} & \frac{1709191}{22248} & \frac{504937}{7416} & \frac{2826673}{44496} & \frac{132799}{11124} \\
&&&&&&&\\
-\frac{1022}{2781} & -\frac{171331}{2472} & -\frac{86549}{22248} & -\frac{65155}{22248} & -\frac{893689}{11124} & -\frac{259999}{3708} & -\frac{1478497}{22248} & -\frac{37589}{2781} \\
&&&&&&&\\
-\frac{2507}{1854} & -\frac{338839}{2472} & -\frac{54739}{7416} & -\frac{38129}{7416} & -\frac{583883}{3708} & -\frac{171893}{1236} & -\frac{981095}{7416} & -\frac{45887}{1854} \\
&&&&&&&\\
-\frac{5359}{11124} & -\frac{362153}{4944} & -\frac{125867}{44496} & -\frac{56581}{44496} & -\frac{1917649}{22248} & -\frac{563347}{7416} & -\frac{3173623}{44496} & -\frac{152365}{11124}
\end{array}\right) \]

On vérifie bien que $N$ est nilpotente avec :

$A= D + N$ et 
 \[ N^2= \left(\begin{array}{cccccccc}
0 & 0 & 0 & 0 & 0 & 0 & 0 & 0 \\
0 & 0 & 0 & 0 & 0 & 0 & 0 & 0 \\
0 & 0 & 0 & 0 & 0 & 0 & 0 & 0 \\
0 & 0 & 0 & 0 & 0 & 0 & 0 & 0 \\
0 & 0 & 0 & 0 & 0 & 0 & 0 & 0 \\
0 & 0 & 0 & 0 & 0 & 0 & 0 & 0 \\
0 & 0 & 0 & 0 & 0 & 0 & 0 & 0 \\
0 & 0 & 0 & 0 & 0 & 0 & 0 & 0
\end{array}\right) \]

\subsection{Méthode de Chevalley}

\paragraph{Préliminaires:}

On a vu dans l'énoncé du théorème de Dunford-Schwarz que dans la décomposition $u=d+n$, $d$ et et $n$ sont des polynômes en $u$. Malgré la méconnaissance des valeurs propres de $u$, la méthode effective de Newton-Raphson nous a permis d'obtenir $d$ et $n$. Mais qu'en est-il du polynôme $Q$ qui génère la composante $d$ ? Existe-t-il aussi une méthode effective qui pourrait nous fournir un tel polynôme ?

La réponse à toutes ces questions existe déjà dans une méthode qui a été évoquée par Claude Chevalley en 1951 et qui n'a pas été largement diffusée pour que tout le monde en fasse usage. Un article de Danielle Couty, Jean Esterle et Rachid Zarouf a essayé de mettre en lumière cette méthode et d'inciter à son enseignement dans les universités. 

J'apporte ma petite pierre à l'édifice en détaillant la preuve de cette méthode avec des outils élémentaires d'arithmétique dans l'espoir qu'elle devienne accessible à tout bon preneur.

Rappelons le principe de la décomposition de Dunford avec les projecteurs spectraux :

Soit $u\in\mathcal{L}(E)$ avec $\chi_u = \prod\limits_{i=1}^r (X-\lambda_i)^{\alpha_i}$  avec pour tout $1\leq i,j \leq r$, $\lambda_i \in \mathbb{K}$  , $\alpha_i\geq1$ et $\lambda_i \neq \lambda_j$ pour $i\neq j$. on pose $N_i=\ker (u-\lambda_i\Id)^{\alpha_i}$.

 D'après le lemme des noyaux et le théorème de Cayley-Hamilton on sait que : 

$E=\bigoplus\limits_{i=1}^p N_i$ et pour tout $i$, la projection spectrale $\pi_i$ sur $N_i$ parallèlement à $\bigoplus\limits_{j\neq i} N_j$ est un polynôme en $u$.

Par le théorème de Bezout on obtient la décomposition de Dunford  $u=d+n$ avec $d=\sum\limits_{i=1}^{r} \lambda_i \pi_i$ où $\pi_i=U_iQ_i(u)$ , $Q_i=\prod\limits_{j=1,j\neq i}^{r} (X-\lambda_j)^{\alpha_j}$ et les $U_i$ sont obtenus grâce à l'Identité de Bezout $\sum\limits_{i=1}^{r} U_iQ_i = 1$.

Donc en désignant par $Q=\sum\limits_{i=1}^{r} \lambda_iU_iQ_i$, on voit bien que $Q$ vérifie le système de congruences : 
$$ (\mathcal{S}) : Q \equiv \lambda_i  \mod  (X-\lambda_i)^{\alpha_i} , 1\leq i \leq r$$.

Puisque les $(X-\lambda_i)^{\alpha_i}$ sont premiers entre eux deux à deux, d'après le lemme chinois,

( $\mathcal{S}$) admet une solution unique $\mod \chi_u$.

On considère $P=\displaystyle\frac{\chi_u}{\chi_u \wedge \chi'_u}=\prod\limits_{i=1}^{r} (X-\lambda_i)$ qu'on arrive à  calculer avec l'algorithme d'Euclide.

Puisque $P$ et $P'$ sont premiers entre eux, il existe alors d'après le théorème de Bezout $U,V \in \mathbb{K}[X]$ qui vérifient $UP + VP' = 1$.

Enfin, la méthode de Chevalley consiste à construire une suite de polynômes qui converge à partir d'un certain rang vers un polynôme solution du système de congruences $ (\mathcal{S})$. Si $Q$ désigne une solution de $ (\mathcal{S})$, on aura dans ce cas $d=Q(u)$ et $n=u-Q(u)$.

\begin{theo}[Chevalley]

On définit une suite $( Q_m)_{m\geq 0}\in\mathbb{K}[X]$ donnée par :
\\
\[ \left\{ \begin{array}{c}
Q_0=X - PV        \\
                          
Q_{m+1}=Q_0 \circ r(Q_m)         
\end{array}\right. \]
\\
$r(Q_m)$ est le reste de la division euclidienne de $Q_m$ par $\chi_u$.

Alors la suite $(Q_m)$ est stationnaire à partir d'une certain rang vers la solution du système~($\mathcal{S}$).

\end{theo}

\paragraph{Preuve:}

On définit \[ \left\{ \begin{array}{cccc}
                               s: & \mathbb{K}[X] & \to &    \mathbb{K}[X]   \\
                          
                                   & f &          \mapsto & f \circ (X-PV)
                               \end{array}\right. \]

\bigskip

Montrons par récurrence que pour tout $k\in\mathbb{N}^*$ : $s^k(f(X)) = f(s^k(X))$.

Pour $k=1$ on a d'une part $s(f(X))=f(X-P(X)V(X))$ et d'autre part $f(s(X)) = f(X \circ (X-P(X)V(X))=f(X-P(X)V(X))$. Donc la propriété est vraie au rang 1.

On suppose que pour $k\in\mathbb{N}^*$ , $s^k(f(X))=f(s^k(X))$ et on a dans ce cas $s^{k+1}(f(X))=s(s^k(f(X)))=s(f(s^k(X)))=f(s^k(s(X)))=f(s^{k+1}(X))$.

\bigskip

Par le développement de Taylor il existe  un élément   $T$  de $\mathbb{K}[X]$ vérifiant :

$s(f(X))=f(X-P(X)V(X))=f(X)-f'(X)P(X)V(X)+P^2(X)V^2(X)T(X)$.

 Par suite, $s(f(X)) \equiv f(X)-f'(X)P(X)V(X) \ \ (\mod \ \  P^2V^2)$.

\bigskip

Donc  \[ \left\{ \begin{array}{ccc}
s(P)(X) \equiv  & P(X)-P'(X)P(X)V(X) &  \mod  P^2V^2    \\
s(P)(X) \equiv  & P(X) (1-P'(X)V(X)) &  \mod  P^2V^2    \\
s(P)(X) \equiv  & P(X)U(X)P(X) & \mod  P^2V^2    \\
s(P)(X) \equiv  &U(X)P^2(X) &  \mod P^2V^2    \\
s(P)(X) \equiv  &0& \mod  P^2  \\ 
                       \end{array}\right. \]

\bigskip

Montrons par récurrence que pour tout $k\in\mathbb{N}^*$ on a $s^k(P) \equiv 0 \ \ (\mod \ \ P^{2^k})$.
\\
En effet pour $k=1$, la propriété vient d'être établie ce qui signifie aussi qu'il existe $A_1\in\mathbb{K}[X]$ tel que $s(P)(X)=P^2(X)A_1(X)$. Autrement dit, $P(X-PV)=P^2(X)A_1(X)$.
\\
Supposons que la propriété est vraie jusqu'à un rang $k\geq 1$, c'est-à-dire qu'il existe $A_k\in\mathbb{K}[X]$ tel que $s^k(P)(X)=P^{2^k}(X)A_k(X)$.

$s^{k+1}(P)(X)=s(s^k(P))(X)=s(P^{2^k}(X)A_k(X))=[P^{2^k}(X)A_k(X)]\circ(X-PV)$

Donc $s^{k+1}(P)(X)=P^{2^k}(X-PV)A_k(X-PV)=[P(X-PV)]^{2^k}A_k(X-PV)$

On a alors, $s^{k+1}(P)(X)=(P^2(X)A_1(X))^{2^k}A_k(X-PV)=P^{2^{k+1}}(X)A_k(X-PV)$.

Finalement on obtient $s^{k+1}(P) \equiv 0\ \ (\mod \ \ P^{2^{k+1}})$.

\bigskip

D'autre part, $s(X) = X-P(X)V(X)$ , c'est-à-dire $s(X) \equiv X \ \ (\mod\ \  P)$ et par récurrence immédiate on a aussi $s^k(X) \equiv X \ \ (\mod \ \ P)$.

\bigskip

Pour $i\neq j$ on a $s^k(X) - \lambda_j \equiv \lambda_i - \lambda_j + (X-\lambda_i) \ \ (\mod\ \ P)$.

Il existe alors $R\in\mathbb{K}[X]$ qui vérifie la relation :

$[s^k(X)-\lambda_j]-(X-\lambda_i) = (\lambda_i - \lambda_j) + RP$

$[s^k(X)-\lambda_j]-(X-\lambda_i) = (\lambda_i - \lambda_j) + R\prod\limits_{l=1}^r (X-\lambda_l)$.
\\
Donc $[s^k(X)-\lambda_j]-(X-\lambda_i) [1 +  R\prod\limits_{l=1, l\neq i}^r (X-\lambda_l)]= (\lambda_i - \lambda_j) $.
\\
On en déduit alors avec le théorème de Bezout que $s^k(X)-\lambda_j$ et $(X-\lambda_i)$ sont premiers entre eux pour tout $i\neq j$.
\\
Donc $\prod\limits_{j=1, j\neq i}(s^k(X)-\lambda_j)$ est premier avec $ (X-\lambda_i)^{2^k}$.
\\
Or $P(s^k(X)) = \prod\limits_{j=1}^r (s^k(X)-\lambda_j) = [\prod\limits_{j=1, j\neq i}^r (s^k(X)-\lambda_j)](s^k(X)-\lambda_i)$ .
\\
Puisque $P(s^k(X))=s^k(P(X)) \equiv 0\ \ (\mod \ \ P^{2^k})$, on en déduit que $P^{2^k}$ divise $P(s^k(X))$.
\\
Autrement dit, $\prod\limits_{i=1}^r (X-\lambda_i)^{2^k}$ divise $ [\prod\limits_{j=1, j\neq i}^r (s^k(X)-\lambda_j)](s^k(X)-\lambda_i)$ .
\\
Soit alors $[\prod\limits_{i=1,i\neq j}^r (X-\lambda_i)^{2^k}](X-\lambda_i)^{2^k}$ qui divise  $ [\prod\limits_{j=1, j\neq i}^r (s^k(X)-\lambda_j)](s^k(X)-\lambda_i)$ .
\\
Par le lemme de Gauss et puisque   $\prod\limits_{j=1, j\neq i}^r(s^k(X)-\lambda_j)$ est premier avec $ (X-\lambda_i)^{2^k}$,  on a $(X-\lambda_i)^{2^k}$ qui divise $(s^k(X)-\lambda_i)$.
\\
Finalement on obtient $s^k(X) \equiv \lambda_i \ \ (\mod\ \  (X-\lambda_i)^{2^k})$.

\bigskip

Soit $m_0$ tel que $2^{m_0+1} \geq \max\{\alpha_i; 1\leq i \leq r\}$.

On a dans ce cas pour $m > m_0$ , $Q_m \equiv s^m(X)  \mod \chi_u$ et il existe dans ce cas $R_m \in\mathbb{K}[X]$ tel que $Q_m(X) = s^m(X)+R_m\chi_u$.
\\
Or il existe aussi $T_m \in\mathbb{K}[X]$ tel que $s^m(X)=\lambda_i+T_m(X-\lambda_i)^{2^m}$.
\\
Donc $Q_m=\lambda_i+R_m\chi_u+T_m(X-\lambda_i)^{2^m}$.

Dans ce cas il existe  $S_m \in\mathbb{K}[X]$ vérifiant $Q_m=\lambda_i+(X-\lambda_i)^{\alpha_i)}S_m$.
\\
\\
Enfin on a  $Q_m \equiv \lambda_i \mod  (X-\lambda_i)^{\alpha_i}$,  c'est-à-dire que $Q_m$ est  solution du système de congruences ($\mathcal{S}$). 
\\
\\
On sait que ($\mathcal{S}$) admet une solution unique $(\mod  \chi_u)$ grâce au lemme chinois et donc $Q_m$ est stationnaire à partir de $m_0$.

\newpage

\paragraph{L'alogorithme en Xcas:} 

On applique la méthode ci-dessus avec Xcas où $A$  la matrice de départ et on obtient à la sortie le polynôme $Q$ vérifiant $D=Q(A)$ ainsi que $D$ et $N$.

\bigskip

\parbox{12cm}{chinois\_effective(A):= \\
\{ local chi,P,Q,u,v,bez,q0,qn,qnn,rqn,rqnn,D,N; \\
  chi:=pcar(A,x);  \\
  P:=simplify(chi/gcd(chi,diff(chi,x)));  \\
  Q:=diff(P,x);  \\
  gcdex(P,Q,x,'u','v');  \\
  bez:=P*u+Q*v;  \\
  simplify(bez);  \\
  v:=simplify(v/bez);  \\
  q0:=simplify(x-P*v);  \\
  qn:=q0;  \\
  rqn:=simplify(rem(qn,chi,x));  \\
  qnn:=simplify(subs(q0,x=rqn));  \\
  rqnn:=simplify(rem(qnn,chi,x));  \\
  while(rqn{\tt\symbol{60}}{\tt\symbol{62}}rqnn)\{ \\
      rqn:=rqnn;  \\
      qn:=simplify(subs(q0,x=rqn));  \\
      rqnn:=simplify(rem(qn,chi,x));  \\
    \};; ;  \\
  D:=simplify(subs(rqn,x=A));  \\
  N:=A-D;  \\
  return(rqn,D,N);  \\
\} }

\paragraph{Exemple :}

On reprend la matrice $A\in\mathcal{M}_8(\mathbb{R})$ qui a été donnée dans l'exemple d'illustration pour la méthode de Newton-Raphson. On applique la méthode de Chevalley avec les restes chinois pour trouver un polynôme $Q$ vérifiant $D=Q(A)$.
\\
\\
Rappelons que $A=\left(\begin{array}{cccccccc}
1 & 163 & 1 & 0 & 162 & 165 & 163 & -1 \\
0 & 209 & -2 & 0 & 210 & 209 & 208 & 2 \\
0 & -86 & -1 & 1 & -87 & -87 & -87 & 1 \\
1 & -22 & 4 & -1 & -23 & -20 & -20 & -4 \\
1 & 155 & 4 & 0 & 153 & 157 & 155 & -3 \\
0 & -111 & -1 & 0 & -110 & -112 & -110 & 1 \\
-1 & -245 & -1 & 0 & -244 & -246 & -245 & 1 \\
-1 & -160 & -3 & 1 & -158 & -162 & -161 & 4
\end{array}\right) \in\mathcal{M}_8(\mathbb{R})$

\bigskip

Avec Xcas on exécute chinois\_effective(A), et on obtient le polynôme :

$$\boxed{Q(x) = \displaystyle\frac{-437 x^{7}+4374 x^{6}-10290 x^{5}+2380 x^{4}-25941 x^{3}-19770 x^{2}+128556 x-28584}{133488}}$$

Puis on calcule $Q(d)$ avec l'instruction D:=subs(Q,x=A) et on obtient le résultat attendu, à savoir :

$D=\left(\begin{array}{cccccccc}
\frac{3719}{16686} & \frac{315275}{3708} & -\frac{117607}{33372} & -\frac{81239}{33372} & \frac{1219219}{16686} & \frac{471985}{5562} & \frac{2927425}{33372} & -\frac{227905}{16686} \\
&&&&&&&\\
-\frac{23185}{16686} & \frac{240737}{3708} & -\frac{357343}{33372} & -\frac{239051}{33372} & \frac{379130}{8343} & \frac{360811}{5562} & \frac{2330431}{33372} & -\frac{208553}{8343} \\
&&&&&&&\\
\frac{6833}{11124} & -\frac{75491}{1648} & \frac{8389}{44496} & \frac{70883}{44496} & -\frac{893137}{22248} & -\frac{338575}{7416} & -\frac{2108287}{44496} & \frac{91001}{11124} \\
&&&&&&&\\
\frac{43421}{33372} & \frac{421813}{14832} & \frac{1135663}{133488} & \frac{451145}{133488} & \frac{2223845}{66744} & \frac{638927}{22248} & \frac{3620027}{133488} & \frac{200327}{33372} \\
&&&&&&&\\
\frac{7133}{11124} & \frac{443321}{4944} & \frac{59371}{44496} & -\frac{45403}{44496} & \frac{1694753}{22248} & \frac{659375}{7416} & \frac{4070207}{44496} & -\frac{166171}{11124} \\
&&&&&&&\\
\frac{1022}{2781} & -\frac{103061}{2472} & \frac{64301}{22248} & \frac{65155}{22248} & -\frac{329951}{11124} & -\frac{155297}{3708} & -\frac{968783}{22248} & \frac{40370}{2781} \\
&&&&&&&\\
\frac{653}{1854} & -\frac{266801}{2472} & \frac{47323}{7416} & \frac{38129}{7416} & -\frac{320869}{3708} & -\frac{132163}{1236} & -\frac{835825}{7416} & \frac{47741}{1854} \\
&&&&&&&\\
-\frac{5765}{11124} & -\frac{428887}{4944} & -\frac{7621}{44496} & \frac{101077}{44496} & -\frac{1597535}{22248} & -\frac{638045}{7416} & -\frac{3990233}{44496} & \frac{196861}{11124}
\end{array}\right) $

\bigskip

$N=\left(\begin{array}{cccccccc}
\frac{12967}{16686} & \frac{289129}{3708} & \frac{150979}{33372} & \frac{81239}{33372} & \frac{1483913}{16686} & \frac{445745}{5562} & \frac{2512211}{33372} & \frac{211219}{16686} \\
&&&&&&&\\
\frac{23185}{16686} & \frac{534235}{3708} & \frac{290599}{33372} & \frac{239051}{33372} & \frac{1372900}{8343} & \frac{801647}{5562} & \frac{4610945}{33372} & \frac{225239}{8343} \\
&&&&&&&\\
-\frac{6833}{11124} & -\frac{66237}{1648} & -\frac{52885}{44496} & -\frac{26387}{44496} & -\frac{1042439}{22248} & -\frac{306617}{7416} & -\frac{1762865}{44496} & -\frac{79877}{11124} \\
&&&&&&&\\
-\frac{10049}{33372} & -\frac{748117}{14832} & -\frac{601711}{133488} & -\frac{584633}{133488} & -\frac{3758957}{66744} & -\frac{1083887}{22248} & -\frac{6289787}{133488} & -\frac{333815}{33372} \\
&&&&&&&\\
\frac{3991}{11124} & \frac{322999}{4944} & \frac{118613}{44496} & \frac{45403}{44496} & \frac{1709191}{22248} & \frac{504937}{7416} & \frac{2826673}{44496} & \frac{132799}{11124} \\
&&&&&&&\\
-\frac{1022}{2781} & -\frac{171331}{2472} & -\frac{86549}{22248} & -\frac{65155}{22248} & -\frac{893689}{11124} & -\frac{259999}{3708} & -\frac{1478497}{22248} & -\frac{37589}{2781} \\
&&&&&&&\\
-\frac{2507}{1854} & -\frac{338839}{2472} & -\frac{54739}{7416} & -\frac{38129}{7416} & -\frac{583883}{3708} & -\frac{171893}{1236} & -\frac{981095}{7416} & -\frac{45887}{1854} \\
&&&&&&&\\
-\frac{5359}{11124} & -\frac{362153}{4944} & -\frac{125867}{44496} & -\frac{56581}{44496} & -\frac{1917649}{22248} & -\frac{563347}{7416} & -\frac{3173623}{44496} & -\frac{152365}{11124}
\end{array}\right) $

\subsection{\'Etude du cas réel}

\subsubsection{Endomorphismes semi-simples}

\begin{defin}[]

On dit que $u\in\mathcal{L}(E)$ est semi-simple si pour tout sous-espace vectoriel $F$ de $E$ stable par $u$, il existe un supplémentaire $S$ de $F$ stable par $u$.

Une matrice $M\in\mathcal{M}_n(\mathbb{K})$ est dite semi-simple si l'endomorphisme $u$ de $\mathbb{K}^n$ dont $M$ est la matrice dans la base canonique de $\mathbb{K}^n$ est semi-simple.

\end{defin}

Pour $u\in\mathcal{L}(E)$, soit $\mu_u =P_1^{\alpha_1} \cdots P_r^{\alpha_r}$ la décomposition du polynôme minimal de $u$ en facteurs irréductibles de $\mathbb{K}[X]$.

\begin{lemme}

Soit $F$ un sous-espace vectoriel de $E$ stable par $u$. Alors :

$$F=\bigoplus\limits_{i=1}^r [\ker P_i^{\alpha_i}(u) \cap F]$$

\end{lemme}

\paragraph{Preuve :}

Notons d'abord $F_i=\ker P_i^{\alpha_i}(u)$. On sait d'après le lemme des noyaux et le théorème de Cayley-Hamilton que $E=F_1\bigoplus \cdots \bigoplus F_r$.

Rappelons que pour tout $i \in $[\![$  1,r$]\!]$ $, $\pi_i$ qui est la projection sur $F_i$ parallèlement à $\bigoplus\limits_{j\neq i} F_j$ est un polynôme en $u$.

Or $F$ est stable par $u$, donc $F$ est aussi stable par $\pi_i$ et dans ce cas on a $\pi_i(F) \subset F$.

Puisque $F \subset E$, $\pi_i(F) \subset \pi_i(E) = F_i$. Alors $\pi_i(F) \subset F\cap F_i$.

D'autre part on sait que $\Id_E = \pi_1+\ldots+ \pi_r$, donc $F \subset \pi_1(F)+\ldots+\pi_r(F)=\pi_1(F)\bigoplus \cdots \bigoplus\pi_r(F)$. 
\\
Soit alors $F \subset (F_1\cap F)\bigoplus \cdots \bigoplus (F_r \cap F)$.
\\
Réciproquement, pour  tout $i \in $[\![$  1,r$]\!]$ $ , $F_i \cap F \subset F$, donc $ (F_1\cap F)\bigoplus \cdots \bigoplus (F_r\cap F) \subset F$.
\\
Avec la double inclusion, on obtient $F=\bigoplus\limits_{i=1}^r (F_i \cap F)$. $\blacksquare$

\begin{lemme}

Si $\mu_u$ est irréductible alors $u$ est semi-simple.

\end{lemme}

\paragraph{Preuve :}

On suppose que $\mu_u$ est irréductible dans $\mathbb{K}[X]$.

Soit $F$ un sous-espace vectoriel de $E$ stable par $u$. Montrons que dans ce cas, il existe un supplémentaire $S$ de $F$ dans $E$ stable par $u$.

Si $F=E$ alors $S=\{0\}$ convient.

Sinon, soit $x_1\in E\backslash F$ et considérons $E_{x_1} = \{P(u)(x_1) ; P\in\mathbb{K}[X]\}$. Soit $y\in E_{x_1}$, il existe alors $P\in\mathbb{K}[X]$ tel que $y=P(u)(x_1)$ et dans ce cas $u(y)=u\circ P(u)(x_1)=Q(u)(x_1)\in E_{x_1}$ où $Q(x)=x P(x)$.

Donc $E_{x_1}$ est stable par $u$. Montrons alors que $E_{x_1} \cap F=\{0\}$.

Pour cela onconsidère $I_{x_1}=\{P\in\mathbb{K}[X]; P(u)(x_1)=0\}$ qui est un idéal de $\mathbb{K}[X]$.

Or $\mu_u(u)(x_1)=0$, donc $\mu_u\in I_{x_1}$ et par suite $I_{x_1}\neq \{0\}$. Il existe alors $\mu_{x_1}$ polynôme unitaire tel que $I_{x_1} =\mu_{x_1}\mathbb{K}[X]$ ( car $\mathbb{K}[X]$ est anneau principal).

Puisque $\mu_u\in I_{x_1}$, $\mu_{x_1}$ divise $\mu_u$. Or $\mu_u$ est irréductible dans $\mathbb{K}[X]$, donc $\mu_{x_1}=\mu_u$ et dans ce cas $\mu_{x_1}$ est irréductible.

Soit $y\in E_{x_1} \cap F$, il existe alors $P\in\mathbb{K}[X]$ tel que $y=P(u)(x_1)$. Supposons que $y\neq 0$, dans ce cas $P\notin I_{x_1}=(\mu_{x_1})$, et donc $\mu_{x_1}$ ne divise pas $P$. Or $\mu_{x_1}$ est irréductible, donc $\mu_{x_1}$ et $P$ sont premiers entre eux.

 D'après le théorème de Bezout, il existe $U,V\in\mathbb{K}[X]$ tels que $UP+V\mu_{x_1}=1$.

Donc $x_1=U(u)\circ P(u)(x_1)+V(u)\circ \mu_{x_1}(u)(x_1)=U(u)\circ P(u)(x_1)=U(u)(y)$.

Or $y\in F$ et $F$ est stable par $u$, donc $x_1=U(u)(y) \in F$, ce qui est absurde car $x_1\in E\backslash F$. On a alors $y=0$ et $E_{x_1} \cap F=\{0\}$.

$E_{x_1}$ et $F$ sont alors en somme directe et $E_{x_1}$ est stable par $u$.

Si $E=F\bigoplus E_{x_1}$ alors $u$ est semi-simple.

Sinon on choisit $x_2\in E\backslash (F\bigoplus E_{x_1})$ et on reprend le même raisonnement en remplaçant $F$ par $F\bigoplus E_{x_1}$.

Puisque $E$ est de dimension finie, au bout d'un nombre fini d'itérations on trouve $x_1, \ldots ,x_k$ tels que $E=F\bigoplus E_{x_1}\bigoplus\cdots\bigoplus E_{x_k}$ où $E_{x_i}$ est stable par $u$ pour tout $1\leq i \leq k$.

Finalement $S= E_{x_1}\bigoplus\cdots\bigoplus E_{x_k}$ est stable par $u$ et vérifie $E=F\bigoplus S$. $\blacksquare$

\begin{prop}

$u$ est semi-simple si et seulement si $\mu_u$ est produit de polynômes irréductibles unitaires distincts deux à deux.

\end{prop}

\paragraph{Preuve :}

\underline{Condition nécessaire :}

Supoosons que $u$ est semi-simple. Soit $\mu_u=P_1^{\alpha_1}\cdots P_r^{\alpha_r}$ la décomposition de $\mu_u$ en facteurs irréductibles unitaires de $\mathbb{K}[X]$. 

Montrons alors que pour tout $i\in $[\![$  1,r$]\!]$ , \alpha_i=1$.

Supposons qu'il existe $i\in $[\![$  1,r$]\!]$ $ tel que $\alpha_i\geq 2$. 
Si $P=P_i$, alors il existe $Q\in\mathbb{K}[X]$ tel que $\mu_u=P^2Q$.

Soit $F=\ker P(u)$. $F$ est stable par $u$ et $u$ est semi-simple, donc il existe un supplémentaire $S$ de $F$ stable par $u$.

Montrons que $PQ(u)$ s'annule sur $S$. Si $x\in S$, alors $PQ(u)(x) \in F$ car $P(u)\bigl(PQ(u)(x)\bigr) = P^2Q(u)(x)=\mu_u(u)(x) =0$ et $PQ(u)(x)\in S$ puisque $S$ est stable par $u$.

Donc $PQ(u)(x) \in F\cap S=\{0\}$ et on a dans ce cas $PQ(u)(x)=0$ et $PQ(u)$ s'annule alors sur $S$.

Si $y\in F=\ker P(u)$, alors $PQ(u)(y)=Q\bigl(P(u)(y)\bigr)=0$. $PQ(u)$ s'annule alors sur $F$. 

Or $E=F\bigoplus S$, donc $PQ(u)$ s'annule sur $E$ ce qui implique que $\mu_u=P^2Q$ divise $PQ$. Ceci est absurde. 
\\
\\
\underline{Condition suffisante :}

Supposons $\mu_u=P_1\cdots P_r$ avec les $P_i$ irréductibles unitaires et distincts deux à deux. Soit $F$ un s.e.v de $E$ stable par $u$.

On note pour tout $i\in $[\![$  1,r$]\!]$, F_i=\ker P_i(u)$. Alors $E=F_1\bigoplus\cdots\bigoplus F_r$ et d'après le lemme 4 , $F=\bigoplus\limits_{i=1}^r\bigl(F_i\cap F\bigr)$.

Pour tout $i\in $[\![$  1,r$]\!]$, F_i$ est stable par $u$. Soit $u_i\in\mathcal{L}(F_i)$ l'endomorphisme induit  de $u$ sur $F_i$. On a $P_i(u_i)=0$ et $P_i$ est irréductible, donc $P_i$ est le polynôme minimal de $u_i$.

D'après le lemme 5 , $u_i$ est semi-simple. Or $F\cap F_i$ est stable par $u_i$, donc il existe un supplémentaire $S_i$ stable par $u_i$ tel que $(F_i\cap F)\bigoplus S_i=F_i$.

On pose alors $S=S_1\bigoplus\cdots\bigoplus S_r$ et on a $E=F_1\bigoplus\cdots\bigoplus F_r= \bigoplus\limits_{i=1}^r \bigl((F_i\cap F)\bigoplus S_i\bigr)=\bigl(\bigoplus\limits_{i=1}^r F_i\cap F\bigr)\bigoplus \bigl(\bigoplus\limits_{i=1}^rS_i\bigr)=F\bigoplus S$.

S est stable par $u$, donc $u$ est semi-simple. $\blacksquare$

\begin{theo}

Soit $M \in \mathcal{M}_n(\mathbb{R})$. 

$M$ est semi-simple dans $\mathcal{M}_n(\mathbb{R})$ si et seulement si $M$ est diagonalisable dans $\mathcal{M}_n(\mathbb{C})$

\end{theo}

\paragraph{Preuve :}

\underline{Condition nécessaire :} Supposons $M$ semi-simple. D'après la proposition précédente $\mu_M=P_1\cdots P_r$ où les $P_i$ sont irréductibles dans $\mathbb{R}[X]$, unitaires et distincts deux à deux.

Montrons que $\mu_M$ n'a que des racines simples dans $\mathbb{C}$. Soit alors $\alpha\in\mathbb{C}$ une racine de $\mu_M$. 
Il existe alors $i\in [\![ 1,r]\!]$ tel que $P_i(\alpha)=0$, soit par exemple $P_1(\alpha)=0$.

Or $P_1$ est irréductible dans $\mathbb{R}[X]$, donc $P_1$ et $P'_1$ sont premiers entre eux et il existe alors $U,V\in\mathbb{R}[X]$ tels que $UP_1+VP'_1=1$.
On alors $UP_1(\alpha)+VP'_1(\alpha)=1$, soit alors $VP'_1(\alpha)=1$ et donc $P'_1(\alpha)\neq 0$. Donc $\alpha$ est une racine simple de $P_1$.

Si $i\neq 1$, $P_1$ et $P_i$ sont premiers entre eux et donc d'après Bezout il existe $U_1,V_1\in\mathbb{K}[X]$ tels que $U_1P_1+V_1P_i=1$, on a alors $U_1P_1(\alpha)+V_1P_i(\alpha)=1$ , soit alors $V_1P_i(\alpha)=1$ et donc $P_i(\alpha)\neq 0$.

Finalement $\alpha$ est racine simple de $\mu_M$. Donc $\mu_M$ est scindé à racines simples dans $\mathbb{C}[X]$ et donc $M$ est diagonalisable dans $\mathcal{M}_n(\mathbb{C})$.
\\
\\
\underline{Condition suffisante :}

Soit $\mu_M=P_1^{\alpha_1}\cdots P_r^{\alpha_r}$ la décomposition de de $\mu_M$ en facteurs irréductibles unitaires dans $\mathbb{R}[X]$.

D'après la proposition précédente, il suffit de montrer que $\alpha_i=1$ pour tout $i\in [\![ 1,r]\!]$.

En effet, $M$ est diagonalisable dans $\mathcal{M}_n(\mathbb{C})$ c'est-à-dire que $\mu_M$ n'a que des racines simples dans $\mathbb{C}$.
 Puisque le polynôme minimal de $M$ dans $\mathbb{C}[X]$ est le même que dans $\mathbb{R}[X]$, $\mu_M$ est produit de polynômes irréductibles unitaires et donc $M$ est semi-simple. $\blacksquare$

\subsubsection{Dunford dans $\mathcal{M}_n(\mathbb{R})$ }

\begin{theo}[Dunford généralisé]

Soit $M \in \mathcal{M}_n(\mathbb{R})$. Alors il existe un unique couple $(S,N) \in \mathcal{M}_n(\mathbb{R})^2$ tel que :

1) $S$ est semi-simple et $N$ nilpotente.

2) $M=S+N$ et $SN=NS$

\end{theo}

\paragraph{Preuve :}

\underline{Unicité :} Supposons qu'un couple $(S,N) \in \mathcal{M}_n(\mathbb{R})^2$ vérifiant 1) et 2) existe. Donc $S$ est diagonalisable sur $\mathbb{C}$ et la décomposition $M=S+N$ est donc celle de Dunford dans $\mathcal{M}_n(\mathbb{C})$. Le couple $(S,N)$ est alors unique.

\underline{Existence :} Il suffit de démontrer que la décomposition de Dunford de $M$ est réelle.

Soit alors $M=D+N$ la décomposition de Dunford dans $\mathcal{M}_n(\mathbb{C})$ avec $D\in\mathcal{M}_n(\mathbb{C})$ diagonalisable et $N\in\mathcal{M}_n(\mathbb{C})$ nilpotente et $DN=ND$.

 Montrons alors que $D,N \in \mathcal{M}_n(\mathbb{R})$.

En effet , $\overline{M}=\overline{D}+\overline{N}$ et $\overline{M}=M$ puisque $M\in\mathcal{M}_n(\mathbb{R})$. On alors $M=\overline{D}+\overline{N}$.

Or $\overline{D}$ est diagonalisable et $\overline{N}$ est nilpotente, donc par unicité, $D=\overline{D}$ et $N=\overline{N}$.

Finalement, $D$, $N\in\mathcal{M}_n(\mathbb{R})$.

\subsubsection{Un algorithme de vérification}

Maintenant qu'on sait que la décomposition d'une matrice $M  \in \mathcal{M}_n(\mathbb{R})$ est de la forme $S+N$ où $S$ est semi-simple dans $\mathbb{R}$, on souhaiterait savoir si $S$ est aussi diagonalisable sur $\mathbb{R}$ sans connaître ses valeurs propres.

On sait que $S$ est toujours diagonalisable sur $\mathbb{C}$ et que $P =\displaystyle\frac{\chi_A}{\chi_A \wedge \chi'_A}$ est un polynôme annulateur de $S$ qui est scindé à racines simples sur $\mathbb{C}$.

Pour que $S$ soit aussi diagonalisable sur $\mathbb{R}$, il suffit que $P$ ait toutes ses racines dans $\mathbb{R}$. Or, le nombre de racines de $P$ est égal à son degré, par suite $S$ est diagonalisable sur $\mathbb{R}$ si le nombre de racines réelles de $P$ est égal au degré de $P$.

Sachant que le théorème de Sturm nous fournit un moyen de compter le nombre de racines réelles d'un polynôme dans $\mathbb{R}[X]$, on va procéder comme suit :

1) On calcule $\chi_A$ puis $ P =\displaystyle\frac{\chi_A}{\chi_A \wedge \chi'_A}$

2) On construit la suite de Sturm $P, P',.......,P_m$ où $P_m=P \wedge P'$.

3) Si $P=\sum\limits_{i=0}^k a_k X^k$,  on calcule la borne de Cauchy $M=1+\max\limits_{0\leq i \leq k-1}\displaystyle\frac{\lvert a_i \rvert}{\lvert a_k \rvert}$

4) On calcule $r=V(-M) - V(M)$ qui est le nombre de racines réelles de $P$.

5) Si $r=\deg P$ alors $S$ est diagonalisable sur $\mathbb{R}$, sinon $S$ est uniquement diagonalisable sur $\mathbb{C}$.

Donc on dispose d'une méthode effective de test de diagonalisabilité sur $\mathbb{R}$ et on a la propriété suivante : 

\begin{prop}
Soit $A\in \mathcal{M}_n(\mathbb{R})$ de polynôme caractéristique $\chi_A$ et $S+N$ sa décomposition de Dunford généralisée . Soit $P=\displaystyle\frac{\chi_A}{\chi_A \wedge \chi'_A}$.

Si $P(X)=X^r+\sum\limits_{k=r-1}^0 a_kX^k$ et $M=1+\max\limits_{0\leq i \leq r-1} a_i $ alors :

$V(-M)-V(M)=r \iff \text{S est diagonalisable sur } \  \mathbb{R}$
\end{prop}

Voici un algorithme écrit pour Xcas qui appelle le programme suite\_sturm( ) déjà détaillé dans la première partie :
\bigskip

\parbox{12cm}{semi\_diag(A):=\{ \\

\{ local chi,P,l,M,d,r; \\

  chi:=simplify(pcar(A,x));  \\

  P:=simplify(chi/gcd(chi,diff(chi,x)));  \\

  l:=convert(symb2poly(P),list);  \\

  M:=1+max(abs(l));  \\

  d:=degree(P);  \\

  r:=suite\_sturm(P,-M,M);  \\

  if (r{\tt\symbol{60}}{\tt\symbol{62}}d) print("A est non diagonalisable sur R");; else print("A est  diagonalisable sur R");; ;  \\
\} }







\end{document}